\documentclass[a4paper,11pt]{book}
\usepackage{amsmath, amsfonts, amssymb}
\usepackage[english, russian]{babel}
\newtheorem{lemma}{Lemma}
\newtheorem{theorem}{Theorem}
\newtheorem{definition}{Definition}

\begin{document}

\begin{center}
\textbf{\large Решение задач Дирихле и Хольмгрена для  трехмерного
сингулярного эллиптического уравнения  методом потенциалов}
\end{center}

\medskip

\begin{center}
\textbf{Т. Г.\,~Эргашев}
\end{center}

\medskip

\begin{quotation}
\textbf{Аннотация}. Потенциалы играют важную роль при решении
краевых задач для эллиптических уравнений. При этом решение ищется
в виде потенциала определенного слоя с неизвестной плотностью, для
определения которой применяется теория интегральных уравнений
Фредгольма второго рода. В свою очередь, такой потенциал
выписывается через фун\-да\-мен\-таль\-ное решение данного
эллиптического уравнения. Используя свойства соответствующих
фундаментальных решений, в прошлом веке была построена теория
потенциала для двумерных эллиптических уравнений с одной и двумя
линиями вырождения. В настоящее время фундаментальные решения
многомерного эллиптического уравнения с несколькими сингулярными
коэффициентами известны. В данной работе мы исследуем потенциалы
двойного и простого слоев для трехмерного эллиптического уравнения
с одним сингулярным коэффициентом и полученные результаты применим
к решению краевых задач Дирихле и Хольмгрена.
\end{quotation}

\noindent \textbf{Ключевые слова.} {Трехмерное эллиптическое
уравнение с одним сингулярным коэффициентом, фундаментальные
решения, теория потенциала, задача Дирихле, задача Хольмгрена.}


\section{Введение} \label{S1}

Многочисленные приложения потенциалов двойного и простого слоев, а
также объемного потенциала можно найти в механике жидкости,
элас\-то\-ди\-на\-ми\-ке, электромагнетизме и акустике, поэтому
теория потенциала за\-ни\-ма\-ет важное место при решении краевых
задач для эллиптических урав\-не\-ний, с помощью которой краевые
задачи удаётся свести к решению интегральных уравнений
\cite{{Mixlin}, {Gunter}, {Kondrat}}.

Различные интересные задачи для уравнения с одним сингулярным
коэффициентом
\begin{equation}
\label{eq0}  u_{xx}+u_{yy} + {\frac{{2\alpha}} {{x}} }u_{x} =
0,\,\,0<2\alpha <1, \,\,\,x>0
\end{equation}
в полуплоскости изучены многими авторами
\cite{{Gel},{Frankl},{Pulkin},{Sm}}. В работах
\cite{{SHC},{Compl},{Ufa},{Tomsk}} исследуются
 потенциалы двойного слоя для обобщенного двуосесимметрического
эллиптического уравнения
\begin{equation}
\label{eq00}  u_{xx}+u_{yy} + {\frac{{2\alpha}} {{x}} }u_{x}+
{\frac{{2\beta}} {{y}} }u_{y} = 0,\,\,0<2\alpha,\,\,2\beta<1
\end{equation}
в области, ограниченной в первой четверти плоскости  $xOy$.

При построении теории потенциала важную роль играют
фундаментальные решения данного эллиптического уравнения. В работе
\cite{MavGar} найдены в яв\-ном виде фундаментальные решения
многомерного сингулярного уравнения Гельмгольца
\begin{equation}
\label{Lob2020} \sum\limits_{i=1}^mu_{x_ix_i} + {\frac{{2\alpha}}
{{x_1}} }u_{x_1} +\lambda u = 0,
\end{equation}
где $m$ -- размерность пространства; \, $\alpha$ и $\lambda$ --
действительные числа, причем $0<2\alpha<1$ и
$-\infty<\lambda<\infty$, а в недавной работе \cite{LobErg} автору
удалось выписать  в явном виде фундаментальные решения
многомерного уравнения Гельмгольца с несколькими сингулярными
коэффициентами. Несмотря на то, что известны фундаментальные
решения даже для многомерного (более двумерного) уравнения, как
уравнение (\ref{Lob2020}), построение теории потенциала
ограничивалось двумерными уравнениями (\ref{eq0})--(\ref{eq00}) и
сравнительно мало работ посвящены когда размерность уравнения
превышает два. Отметим лишь работы \cite{{Mav},{Muh}}.

В настоящей работе мы рассмотрим уравнение
\begin{equation}
\label{eq1} E(u) \equiv u_{xx}+u_{yy}+u_{zz} + {\frac{{2\alpha}}
{{x}} }u_{x}  = 0\,\,  \left(0 < 2\alpha < 1\right)
\end{equation}
в полупространстве $x > 0$.  Для уравнения (\ref{eq1}) построим
теорию потенциала и применим ее к решению основных краевых задач.

При необходимости будем использовать запись $R_3^+=\left\{(x,y,z):
x>0\right\}$ для обозначения полупространства $x>0$. Прежде чем
перейти к изложению ос\-нов\-ных результатов приведем необходимые
сведения о специальных функциях.

Пусть $t\in \mathbb{C}$. Гамма-функция $\Gamma(t)$ определяется в
виде интеграла Эйлера второго рода \cite[гл. 1, \S 1.1, формула
(1)]{Bateman}
\begin{equation}
\label{Gam} \Gamma(t)=\int_0^\infty s^{t-1}e^{-s}ds,
\end{equation}
который сходится при всех $t\in \mathbb{C}$, для которых  $\Re
t>0$.

Интегрирование по частям выражения (\ref{Gam}) приводит к
рекуррентной формуле
\begin{equation}
\label{Recurr} \Gamma(t+1)=t\Gamma(t).
\end{equation}
Переписав формулу (\ref{Recurr}) в виде
\begin{equation}
\label{Recurr1} \Gamma(t-1)=\frac{\Gamma(t)}{t-1},
\end{equation}
мы получим выражение, позволяющее определить гамма-функцию от
отрицательных аргументов, для которых определение (\ref{Gam})
неприемлемо. Формула (\ref{Recurr1}) показывает, что $\Gamma(t)$
имеет в точках $t=0,\,-1,\,-2,\,-3,...$ простые полюсы.

Символ Похгаммера $(t)_n$ при целых $n$ определяется равенством
\begin{equation}
\label{Pox} (t)_n=t(t+1)...(t+n-1),\,\, n=1,2,...; \,\,(t)_0\equiv
1.
\end{equation}

Справедливы равенства $(t)_n=(-1)^n(1-n-t)_n,\,\,(1)_n=n! $   и
\begin{equation}
\label{GamPox} (t)_n=\frac{\Gamma(t+n)}{\Gamma(t)}
\end{equation}

Равенство (\ref{GamPox}) можно использовать для введения символа
$(t)_n$ при действительных (комплексных) $n$.

Гипергеометрическая функция Гаусса  определяется внутри круга
$|t|<1$ как сумма гипергеометрического ряда \cite[гл.2, \S 2.1,
формула (2)]{Bateman}
\begin{equation} \label{Gauss}
F\left( {a,b;c;t} \right) = {\sum\limits_{k = 0}^{\infty}
{{\frac{{(a)_{k} (b)_{k}}} {{k!(c)_{k}}} }t^{k}}}  ,
\end{equation}
а при $|t|\geq 1$ получается аналитическим продолжением этого ряда
. В формуле (\ref{Gauss}) параметры $a,\,b,\,c$  и  переменная $t$
могут быть комплексными, причем $c\neq 0,\,-1,\,...,$  а $(a)_k$
есть символ Похгаммера  (\ref{Pox}).

Фундаментальные решения уравнения (\ref{eq1}) имеют вид
\cite{{MavGar},{LobErg},{Olev},{SH}}:
\begin{equation}
\label{eq2} q_{1} \left( {x,y,z;\xi,\eta,\zeta}  \right)
=\frac{1}{2\pi} r^{-2\alpha-1} F\left( {\alpha +
{\frac{{1}}{{2}}},\alpha ;2\alpha ;\sigma} \right),
\end{equation}
\begin{equation}
\label{eq3} q_{2} \left( {x,y,z;\xi,\eta,\zeta}  \right)
=\frac{1}{2\pi} r^{2\alpha-3}x^{1 - 2\alpha} \xi^{1 - 2\alpha}
F\left( {{\frac{{3}}{{2}}} - \alpha ,1 - \alpha ;2 - 2\alpha
;\sigma} \right),
\end{equation}
где
\begin{equation*} \label{eq41}
\sigma = 1 - {\frac{{r_{1}^{2}}} {{r^{2}}}};\,\,\, \left.
{{\begin{array}{*{20}c}
 { r_1^2} \hfill \\
 { r^2} \hfill \\
\end{array}}}  \right\}={{\left( {x\pm
\xi}  \right)^{2}}}+{{\left( {y - \eta}  \right)^{2}}}+{{\left( {z
- \zeta}  \right)^{2}}}.
\end{equation*}
Эти функции по переменным $(x,y,z)$ являются решениями уравнения
(\ref{eq1}), имеют особенность порядка $\displaystyle\frac{1}{r}$
при $r \to 0$ и, следовательно, действительно являются
фундаментальными решениями уравнения (\ref{eq1}). Нетрудно видеть,
что

\begin{equation}
\label{eq5} {\left. \left({{x^{2\alpha}\frac{{\partial q_{1}\left(
{x,y,z;\xi,\eta,\zeta}  \right)}}{{\partial x }}}}\right)
\right|}_{x = 0} ={\left. \left({{\xi^{2\alpha}\frac{{\partial
q_{1}\left( {x,y,z;\xi,\eta,\zeta}  \right)}}{{\partial \xi
}}}}\right) \right|}_{\xi = 0}= 0,
\end{equation}

\begin{equation} \label{eq6}
{\left. {q_{2} \left( {x,y,z;\xi,\eta,\zeta}  \right)} \right|}_{x
= 0} = {\left. {q_{2} \left( {x,y,z;\xi,\eta,\zeta}  \right)}
\right|}_{\xi = 0}=0
\end{equation}
для всех $y$, $z$, $\eta$ и $\zeta$.

\section{Формулы Грина} \label{S2}

Рассмотрим тождество
\begin{equation*}
\label{eq7} x^{2\alpha}  \left[ uE(v) - vE(u) \right] =
\end{equation*}
\begin{equation}
\label{eq7} = \frac{\partial} {\partial x}\left[{x^{2\alpha}\left(
v_{x}  u - vu_{x} \right)}\right]+x^{2\alpha}\frac{\partial}
{\partial y}{\left( v_{y}  u - vu_{y}
\right)}+x^{2\alpha}\frac{\partial} {\partial z}{\left( v_{z}  u -
vu_{z} \right)} .
\end{equation}

Интегрируя обе части последнего тождества по области $D$,
расположенной в полупространстве $x>0$, и пользуясь формулой
Гаусса-Остроградского, по\-лу\-чим

\begin{equation}
\label{eq8} \int\int\int_D  x^{2\alpha} \left[ uE(v)-vE(u)
\right]dxdydz = \int\int_S
{{\left({uB_n^{\alpha}[v]-vB_n^{\alpha}[u]}\right)}}dS,
\end{equation}
где $S\,$ граница области $D$, $n$ -- внешняя нормаль к
поверхности $S$ и $$
B_n^{\alpha}[\,\,\,]=x^{2\alpha}\left(cos(n,x)\cdot\frac{\partial}{\partial
x}+cos(n,y)\cdot\frac{\partial}{\partial y}+
cos(n,z)\cdot\frac{\partial}{\partial z}\right) $$ -- конормальная
производная.

Формула Грина (\ref{eq8}) выводится при следующих предположениях:
функции $u(x,y,z)$, $v(x,y,z)$  и их частные производные первого
порядка непрерывны в замкнутой области $\overline{D}$, частные
производные второго порядка непрерывны внутри $D$ и интегралы по
$D$, содержащие $E(u)$ и $E(v)$, имеют смысл. Если $E(u)$ и $E(v)$
не обладают непрерывностью вплоть до $S$, то это -- несобственные
интегралы, которые получаются как пределы по любой
последовательности областей $D_n$, которые содержатся внутри $D$ ,
когда эти области $D_n$ стремятся к $D$, так что всякая точка,
находящаяся внутри $D$, попадает внутрь областей $D_n$, начиная с
некоторого номера $n$.

Если $u$ и $v$ суть решения уравнения (\ref{eq1}), то из формулы
(\ref{eq8}) имеем
\begin{equation}
\label{eq9} {\int\int_{S}
{{\left({uB_n^{\alpha}[v]-vB_n^{\alpha}[u]}\right)}}dS}  = 0.
\end{equation}

Полагая в формуле (\ref{eq8})  $v\equiv 1$ и заменяя  $u$ на
$u^2$, получим
\begin{equation}
\label{eq10} {\int\int\int_{D}  {x^{2\alpha} \left[{{\left(
{{\frac{{\partial u}}{{\partial x}} }} \right)}}^{2}+{{\left(
{{\frac{{\partial u}}{{\partial y}} }} \right)}}^{2}+{{\left(
{{\frac{{\partial u}}{{\partial z}} }}
\right)}}^{2}\right]dxdydz}}  = {\int\int_{S}uB_n^{\alpha}[u]dS} ,
\end{equation}
где $u\left( {x,y,z} \right)$ -- решение уравнения (\ref{eq1}).

Наконец, из формулы (\ref{eq8}), полагая  $v\equiv1$ , будем иметь
\begin{equation}
\label{eq11} {\int\int_{S} B_n^{\alpha}[u]dS}  = 0,
\end{equation}
т.е интеграл от конормальной производной решения уравнения
(\ref{eq1}) по зам\-кну\-той поверхности $S$ области $D$ равен
нулю.

\newpage

\section{Потенциал двойного слоя} \label{S3}

Пусть $\Gamma$ -- поверхность Ляпунова \cite[гл.18, \S 1]{Mixlin},
лежащая в полупространстве $x>0$ и пусть $D$ -- область,
ограниченная односвязной открытой областью $X$ плоскости $yOz$ и
поверхностью $\Gamma$. Общую границу плоской области $X$ и
поверхности $\Gamma$ обозначим через $\gamma$.

Параметрическое уравнение поверхности $\Gamma$ пусть будет
$x=x(s,t)$, $y=y(s,t)$, $z=z(s,t),$ $(s,t)\in \overline{\Phi}$,
где $\Phi:=(s_1,\,s_2)\times(t_1,\,t_2)$ -- область изменения $s$
и $t$. Тогда параметрическое уравнение плоской кривой $\gamma$
будет иметь вид $y=y(s,t_0),\,z=z(s,t_0)$, где $t_0\in[t_1,\,t_2]$
-- фиксированное число, удовлетворяющее уравнению $x(s,t_0)=0$ при
любых значениях $s\in [s_1,\,s_2]$. Относительно поверхности
$\Gamma$ будем предполагать, что:

1) функции  $x(s,t)$, $y(s,t)$ и $z(s,t)$ имеют непрерывные
частные производные первого порядка по $s$ и $t$  в $
\overline{\Phi}$, не обращающиеся одновременно в нуль
($x_s^2+y_s^2+z_s^2\neq 0$, $x_t^2+y_t^2+z_t^2\neq 0$);

2) при стремлении точек поверхности $\Gamma$ к точкам  кривой
$\gamma$ поверхность $\Gamma$ образует прямой угол с плоскостью
$x=0$.

Координаты переменной точки на поверхности $\Gamma$ будем
обозначать через $\left(\xi,\eta,\zeta\right)$, где
$\xi=\xi(\theta,\vartheta)$, $\eta=\eta(\theta,\vartheta)$,
$\zeta=\zeta(\theta,\vartheta),$ $(\theta,\vartheta)\in
\overline{\Phi}$.

\subsection{Потенциал двойного слоя}\label{s3.1}
Рассмотрим интеграл
\begin{equation}
\label{eq302} w^{(1)}(x,y,z) = {\int\int_{\Gamma} {\mu _{1} \left(
\theta,\vartheta \right)B_\nu^{\alpha}\left[q_{1} \left(
{\xi,\eta,\zeta;x,y,z} \right)\right]d\theta d\vartheta}}  ,
\end{equation}
где $\mu _{1} \left( \theta,\vartheta \right) \in C\left(
{\overline {\Gamma}}   \right)$ и $q_{1} \left( {\xi ,\eta,\zeta;
x,y,z} \right)$ -- фундаментальное решение уравнения (\ref{eq1}),
определенное формулой (\ref{eq2}). Здесь
\begin{equation} \label{eq303}
B_\nu^{\alpha}[\,\,\,]=\xi^{2\alpha}\left(cos(\nu,\xi)\cdot\frac{\partial}{\partial
\xi}+ cos(\nu,\eta)\cdot\frac{\partial}{\partial
\eta}+cos(\nu,\zeta)\cdot\frac{\partial}{\partial \zeta}\right),
\end{equation}
$\nu$ -- внешняя нормаль к поверхности $\Gamma$.

\begin{definition}\label{de1}
Интеграл (\ref{eq302}) будем называть  \textit{потенциалом
двойного слоя с плотностью} $\mu _{1} \left(
\theta,\vartheta\right)$. \end{definition}

Очевидно, что $w^{({1})} \left( x,y,z \right)$ есть регулярное
решение уравнения (\ref{eq1}) в любой области, лежащей в
полупространстве $x>0$, не имеющей общих точек ни с поверхностью
$\Gamma$, ни с плоскостью $yOz$. Как и в случае логарифмического
потенциала, можно показать существование потенциала двойного слоя
(\ref{eq302}) в точках поверхности $\Gamma$ для ограниченной
плотности $\mu_{1}\left(\theta,\vartheta \right)$.

\begin{lemma} \label{L1}{Справедлива следующая формула:}
     \begin{equation} \label{eq304}
w_{1}^{({1})} (x,y,z) \equiv
{\int\int_{\Gamma}B_\nu^{\alpha}\left[q_{1} \left(
{\xi,\eta,\zeta;x,y,z} \right)\right]d\Gamma} = \left\{
{{\begin{array}{*{20}c}
 { - 1,\,\,\,\,(x,y,z) \in D \cup X,} \hfill \\
 { - {\displaystyle\frac{{1}}{{2}}},\,\,\,\,(x,y,z) \in \Gamma \cup \gamma ,} \hfill \\
 {\,\,\,\,\,\,0,\,\,\,\,\,(x,y,z) \notin \bar {D}.} \hfill \\
\end{array}}}  \right.
\end{equation}
\end{lemma}

 \textbf{Доказательство.} Процесс доказательства состоит из нескольких случаев.

\textbf{1-случай.} Пусть точка $(x,y,z)$ лежит вне области $D$ и
над плоскостью $yOz$. Тогда  $q_{1} \left( {\xi ,\eta,\zeta;
x,y,z} \right)$  есть регулярное решение уравнения (\ref{eq1})
внутри области $D$ с непрерывными производными всех порядков
вплоть до поверхности $\Gamma$ , и в силу (\ref{eq11})
\begin{equation*}
\label{eq12000} w_1^{\left( {1} \right)}\left( {x},y,z \right)
\equiv {\int\int_{\Gamma}  {B_\nu^{\alpha}[q_{1} \left(
{\xi,\eta,\zeta;x,y,z} \right)]d\theta d\vartheta}}=0.
\end{equation*}

\textbf{2-случай.} Пусть точка $(x,y,z)$ находится внутри $D$.
Вырежем из области $D$ шар малого радиуса  $\rho$ с центром в
точке $(x,y,z)$ и обозначим через $D_\rho$ оставшуюся часть
области $D$, а через $C_\rho$ сферу вырезанного шара. В области
$D_\rho$  функция $q_{1} \left( {\xi ,\eta,\zeta; x,y,z}\right)$
-- регулярное решение уравнения (\ref{eq1}) и, согласно
(\ref{eq5}) и (\ref{eq11}), мы имеем
\begin{equation*}
{\int\int_{\Gamma}{{B_\nu^{\alpha}[q_{1} \left(
{\xi,\eta,\zeta;x,y,z} \right)] }d\theta
d\vartheta}}+{\int\int_{C_{\rho}} {  {B_\nu^{\alpha}[q_{1} \left(
{\xi,\eta,\zeta;x,y,z} \right)] }dC_{\rho}} }=0,
\end{equation*}
т.е.
\begin{equation} \label{eq305}
 w_{1}^{(1)} (x,y,z) = -{\int\int_{C_{\rho}} {B_\nu^{\alpha}[q_{1} \left( {\xi,\eta,\zeta;x,y,z} \right)]dC_{\rho}} }.
\end{equation}

Вычислим производную по нормали $\nu$ от фундаментального решения
$q_{1} \left( {\xi,\eta,\zeta; x,y,z} \right)$. Применяя
последовательно формулу для вычисления произ\-вод\-ной
гипергеометрической функции Гаусса \cite[гл.2, \S 2.8, формула
(20)]{Bateman}
$$
\frac{d}{d\sigma}F(a,b;c;\sigma)=\frac{ab}{c}F(a+1,b+1;c+1;\sigma)
$$
и смежное соотношение
$$
\frac{b}{c}\sigma
F(a+1,b+1;c+1;\sigma)=F(a+1,b;c;\sigma)-F(a,b;c;\sigma),
$$
получим
\begin{equation*}
{\frac{{\partial q_{1} \left( {\xi,\eta,\zeta; x,y,z}
\right)}}{{\partial \xi}} } =\frac{1+2\alpha}{2\pi}  (x-\xi)r^{-
2\alpha -3}F\left({\alpha +{\frac{{3}}{{2}}},\alpha ;2\alpha
;\sigma}  \right)
\end{equation*}
\begin{equation} \label{eq13}
-\frac{1+2\alpha}{2\pi} x r^{- 2\alpha -3}F\left( {\alpha +
{\frac{{3}}{{2}}},1 + \alpha ;1 + 2\alpha ;\sigma} \right).
\end{equation}

В результате применения формулы дифференцирования \cite[гл.2, \S
2.8, фор\-му\-ла (21)]{Bateman}
$$
\frac{d}{d\sigma}\left[\sigma^aF(a,b;c;\sigma)\right]=a\sigma^{a-1}F(a+1,b;c;\sigma)
$$
находим
\begin{equation}
\label{eq14} {\frac{{\partial q_{1} \left( {\xi,\eta,\zeta; x,y,z}
\right)}}{{\partial \eta}} } =\frac{1+2\alpha}{2\pi} (y-\eta)r^{-
2\alpha -3}F\left( {\alpha + {\frac{{3}}{{2}}},\alpha ;2\alpha
;\sigma} \right),
\end{equation}

\begin{equation}
\label{eq141}{\frac{{\partial q_{1} \left( {\xi,\eta,\zeta; x,y,z}
\right)}}{{\partial \zeta}} } =\frac{1+2\alpha}{2\pi}(z-\zeta)r^{-
2\alpha -3}F\left( {\alpha + {\frac{{3}}{{2}}},\alpha ;2\alpha
;\sigma} \right).
\end{equation}

Используя (\ref{eq13}), (\ref{eq14}) и (\ref{eq141}),  с учетом
(\ref{eq303}), имеем
\begin{equation*}
B_\nu^{\alpha}[q_{1} \left( {\xi,\eta,\zeta;x,y,z} \right)]
=\frac{1+2\alpha}{2\pi}  r^{-2\alpha-1}F\left( {\alpha +
{\frac{{3}}{{2}}},\alpha ;2\alpha ;\sigma}
\right)B_\nu^{\alpha}\left[ \ln \frac{1}{r} \right]
\end{equation*}
\begin{equation}
\label{eq144} -\frac{1+2\alpha}{2\pi} x r^{- 2\alpha
-3}\xi^{2\alpha}F\left( {\alpha + {\frac{{3}}{{2}}},1 + \alpha ;1
+ 2\alpha ;\sigma}  \right)\cos \left( {\nu,\xi}  \right).
 \end{equation}

Далее, применив известную формулу \cite[гл.2, \S 2.9, формула
(2)]{Bateman}
\begin{equation*}
F\left( {a,b;c;x} \right) = \left( {1 - x} \right)^{ - b}F\left(
{c - a,b;c;{\frac{{x}}{{x - 1}}}} \right)
\end{equation*}
к каждой гипергеометрической функции в (\ref{eq144}), получим
\begin{equation*}
B_\nu^{\alpha}[q_{1} \left( {\xi,\eta,\zeta;x,y,z} \right)] =
-\frac{1+2\alpha}{4\pi r_1^{2\alpha}r}F\left( {\alpha
-{\frac{{3}}{{2}}},\alpha ;2\alpha ;1-\frac{r^2}{r_1^2}}
\right)B_\nu^{\alpha}\left[\ln r^2\right]
\end{equation*}

\begin{equation}
\label{eq306} -\frac{1+2\alpha}{2\pi
rr_1^{2\alpha+2}}x\xi^{2\alpha}F\left( {\alpha -
{\frac{{1}}{{2}}}, 1+\alpha ;1 + 2\alpha ;1-\frac{r^2}{r_1^2}}
\right)\cos\left( {\nu,\xi}  \right).
 \end{equation}

Теперь равенство  (\ref{eq305}) может быть написано так:
\begin{equation}
\label{eq17} w_{1}^{(1)} (x,y,z) = i_{1} (x,y,z) + j_{1} (x,y,z),
\end{equation}
где
\begin{equation}
\label{eq171} i_{1} (x,y,z) = -\frac{1+2\alpha}{2\pi}
{\int\int_{C_{\rho}} \frac{1}{r_1^{2\alpha}r}  F\left( {\alpha
-{\frac{{3}}{{2}}},\alpha ;2\alpha ;1-\frac{r^2}{r_1^2}}
\right)B_\nu^{\alpha}\left[ \ln \frac{1}{r} \right] dC_{\rho} }  ,
\end{equation}

\begin{equation*}
 j_{1} (x,y,z) =\frac{1+2\alpha}{2\pi}
 x{\int\int_{C_{\rho}}  \frac{\xi^{2\alpha}}{r_1^{2\alpha}r}F\left( {\alpha -
{\frac{{1}}{{2}}}, 1+\alpha ;1 + 2\alpha ;1-\frac{r^2}{r_1^2}}
\right)\cos \left( {\nu,\xi} \right)dC_{\rho} } .
\end{equation*}

Преобразуем правую часть равенства (\ref{eq171}). На сфере нормаль
направлена против радиуса. Отсюда
\begin{equation}
\label{eq17100} i_{1} (x,y,z) =  -\frac{1+2\alpha}{2\pi}
{\int\int_{C_{\rho}} \frac{\xi^{2\alpha}}{r_1^{2\alpha}r^2}
F\left( {\alpha -{\frac{{3}}{{2}}},\alpha ;2\alpha
;1-\frac{r^2}{r_1^2}} \right)dC_{\rho} }.
\end{equation}

Вводя сферические координаты
\begin{equation}
\label{eq18} \xi=x+\rho \cos\varphi, \,\,\,\eta=y+\rho \sin\varphi
\cos\psi,\,\,\,\zeta=z+\rho \sin\varphi \sin\psi
\end{equation}
\begin{equation*}
\left(\rho\geq 0, \,\,
0\leq\varphi\leq\pi,\,\,0\leq\psi\leq2\pi\right),
\end{equation*}
в интеграле (\ref{eq17100}), получим
\begin{equation*}
 i_{1} (x,y,z)  = -\frac{1+2\alpha}{2\pi}{\int\limits_{0}^{2\pi}  {d\psi}}
{\int\limits_{0}^{\pi}}\left[\frac{x^2+2x\rho
\cos\varphi+\rho^2\cos^2\varphi}{4x^2+4x\rho\cos\varphi+\rho^2}\right]^{\alpha}
\end{equation*}
\begin{equation}
\label{eq19}
\cdot F\left( {\alpha -{\displaystyle\frac{{3}}{{2}}},\alpha ;2\alpha ;\frac{4x^2+4x\rho\cos\varphi}{4x^2+4x\rho\cos\varphi+\rho^2}}   \right){\sin\varphi}d\varphi. \\
\end{equation}

Теперь в правой части равенства (\ref{eq19}) переходим к пределу
при $\rho\to 0$. Используя формулу суммирования для
гипергеометрической функции Гаусса \cite[гл.2, \S 2.8, формула
(46)]{Bateman}
\begin{equation*}
F(a,b;c;1) = {\frac{{\Gamma (c)\Gamma (c - a - b)}}{{\Gamma (c -
a)\Gamma (c - b)}}}, \,\,Re(c - a - b) > 0;\,\,c \ne 0, - 1, -
2,..., \end{equation*} получим
\begin{equation} \label{eq193}
{\mathop {\lim} \limits_{\rho \to 0}} i_{1} (x,y,z) = -1.
\end{equation}

Еще проще доказывается, что
\begin{equation} \label{eq194}
{\mathop {\lim} \limits_{\rho \to 0}} j_{1} (x,y,z) = 0.
\end{equation}

Подставляя теперь (\ref{eq193}) и  (\ref{eq194}) в (\ref{eq17}),
получим
\begin{equation*}
\label{eq1777} w_{1}^{(1)} (x,y,z) = -1, \,\,\, (x,y,z)\in D.
\end{equation*}

\textbf{3-случай.} Пусть теперь точка $(x,y,z)$ совпадает с
некоторой точкой $M_0\left(x_0,y_0,z_0\right)$, лежащей на
поверхности $\Gamma$. Проведем сферу малого радиуса $\rho$ с
центром в точке $M_0$. Эта сфера вырежет часть $\Gamma_\rho$
поверхности $\Gamma$. Оставшуюся часть поверхности обозначим через
${\Gamma_{\rho}^*}$. Мы имеем
\begin{equation}
\label{eq308} w_{1}^{\left( {1} \right)}
\left({x_0,y_0,z_0}\right) = {\mathop {\lim }\limits_{\rho \to 0}}
{\int\int_{{\Gamma_{\rho}^*}}B_\nu^{\alpha}\left[{{q_{1} \left(
{\xi,\eta,\zeta; x_0,y_0,z_0} \right)}}\right]d{\Gamma_{\rho}^*}}
.
\end{equation}

Обозначим через $C_{\rho}^*$  часть сферы $C_\rho$, лежащей внутри
области $D$ и рассмотрим область, ограниченную поверхностями
$\Gamma_{\rho}^*$, $C_{\rho}^*$ и плоской областью $X$ плоскости
$yOz$. Так как точка $M_0$ лежит вне этой области, то в этой
области $q_1(\xi,\eta,\zeta; x,y,z)$ -- регулярное решение
уравнения (\ref{eq1}), и в силу (\ref{eq11})  мы имеем
\begin{equation*}
{\int\int_{{\Gamma_{\rho}^*}} B_\nu^{\alpha}\left[{{q_{1} \left(
{\xi,\eta,\zeta; x_0,y_0,z_0} \right)}}\right]d{\Gamma_{\rho}^*}}
=
\end{equation*}
\begin{equation}
\label{eq309}  = {\int\int_{{C_{\rho}^*}}
{B_\nu^{\alpha}\left[{{q_{1} \left( {\xi,\eta,\zeta; x_0,y_0,z_0}
\right)}}\right]d{C_{\rho}^*}}} .
\end{equation}

Подставляя (\ref{eq309})  в (\ref{eq308}), получим
\begin{equation*}
\label{eq234} w_{1}^{\left( {1} \right)}
\left({x_0,y_0,z_0}\right) =- {\mathop {\lim }\limits_{\rho \to
0}} {\int\int_{{C_{\rho}^*}} {B_\nu^{\alpha}\left[{{q_{1} \left(
{\xi,\eta,\zeta; x_0,y_0,z_0} \right)}}\right]d{C_{\rho}^*}}} .
\end{equation*}
Вводя снова сферические координаты (\ref{eq18}) с центром в точке
$M_0$, получим
\begin{equation*}
\label{eq301} w_{1}^{(1)} (x,y,z) = -\frac{1}{2}, \,\,\,
(x,y,z)\in\Gamma.
\end{equation*}

\textbf{4-случай.}  Положим, наконец, что точка $(x,y,z)$
находится на плоскости $yOz$. Проведем плоскость $x=\delta$
($\delta>0$ достаточно мало) и рассмотрим область $D_\delta$,
которая есть часть области $D$, лежащая над плоскостью $x=\delta$.
Применяя формулу (\ref{eq11}), получим

\begin{equation*}
\label{eq311} w_{1}^{\left( {1} \right)} \left({0,y,z}\right) =
{\int\int_{{H_{\delta}}} B_\nu^{\alpha}[q_{1} \left(
{\xi,\eta,\zeta;0,y,z} \right)]d{H_{\delta}}} +
\end{equation*}
\begin{equation}
\label{eq311} +{\int\int_{X_{\delta}}}\left.\left[\xi^{2\alpha}
{\frac{{\partial q_{1} \left( {\xi,\eta,\zeta; 0,y,z}
\right)}}{{\partial \xi}}}\right]\right|_{\xi=\delta}dX_\delta ,
\end{equation}
где $H_\delta$ -- часть поверхности $\Gamma$, на\-ход\-я\-ща\-яся
ниже плоскости $x=\delta$, а  $X_\delta$ -- сечение области $D$
плоскостью $x=\delta$, т.е. $X_\delta$ -- плоская область,
ограниченная замкнутой кривой $\gamma_\delta:
y=y\left(s,t_\delta\right)$, $z=z\left(s,t_\delta\right)$, здесь
$t_\delta\in[t_1,\,t_2]$ определяется из урав\-не\-ния
$x\left(s,t_\delta\right)=\delta$ при любых значениях $s\in
[s_1,\,s_2]$.

Нетрудно видеть, что интеграл
$$
{\int\int_{{H_{\delta}}} { B_\nu^{\alpha}[q_{1} \left(
{\xi,\eta,\zeta;0,y,z} \right)]d{H_{\delta}}}}
$$
при $\delta\to 0$ стремится к нулю при любом значении $y$ и $z$.

Пусть $X_\delta=\left\{(y,z): p(\delta)< y < q(\delta),\,\,
h(y,\delta)<z<k(y,\delta)\right\}$,\, где  $p(\delta)$ и
$q(\delta)$ -- некоторые конечные действительные числа, зависящие
от $\delta$, а $z=h(y,\delta)$ и $z=k(y,\delta)$ -- непрерывные
функции на отрезке $[p(\delta),\,q(\delta)]$. Далее, согласно
(\ref{eq13}) выражение (\ref{eq311}) можно записать в виде

$w_{1}^{\left( {1} \right)} \left({0,y,z}\right) =$
\begin{equation}
\label{eq313} =-\frac{1+2\alpha}{2\pi}{\mathop {\lim
}\limits_{\delta\to 0}}\delta^{1+2\alpha}
\int\limits_{p(\delta)}^{q(\delta)}\int\limits_{h(\zeta,\delta)}^{k(\zeta,\delta)}\left[\delta^2+(\eta-y)^2+(\zeta-z)^2\right]^{-\alpha-\frac{3}{2}}d\eta
d\zeta.
\end{equation}
Преобразуем выражение (\ref{eq313}). Вместо $\eta$ и $\zeta$
введем новые переменные ин\-тег\-ри\-ро\-ва\-ния
$t=\displaystyle\frac{\eta-y}{\delta}$   и
$s=\displaystyle\frac{\zeta-z}{\delta}$. Совершая замену
переменных, получим

\begin{equation}
\label{eq314} w_{1}^{\left( {1} \right)} \left({0,y,z}\right) =
-\frac{1+2\alpha}{2\pi}{\mathop {\lim }\limits_{\delta\to 0}}
\int_{\alpha_1}^{\alpha_2}\int_{\beta_1}^{\beta_2}\left(1+t^2+s^2\right)^{-\alpha-\frac{3}{2}}dt
ds,
\end{equation}
где
 $$\alpha_1=\displaystyle\frac{p(\delta)-y}{\delta},\,\,\alpha_2=\displaystyle\frac{q(\delta)-y}{\delta},\,\,
 \beta_1=\displaystyle\frac{h(\eta+t\delta,\delta)-z}{\delta},\,\beta_2=\displaystyle\frac{k(\eta+t\delta,\delta)-z}{\delta}. $$
Известно, что  \cite[гл.4, \S 4.6, формула 4.638.3]{Grad}
\begin{equation}
\label{eq316}
\int\limits_{-\infty}^{\infty}\int\limits_{-\infty}^{\infty}\left(1+t^2+s^2\right)^{-\alpha-\frac{3}{2}}dtds
=\frac{2\pi}{1+2\alpha}.
\end{equation}

Таким образом, из формулы (\ref{eq314}) в силу условий 1) и 2),
наложенных на поверхность $\Gamma$ и формулы (\ref{eq316}) будем
иметь
\begin{equation*}
w_{1}^{\left( {1} \right)} \left( 0,y,z \right) =
 \left\{
{{\begin{array}{*{20}c}
 { - 1,\,\,\,\,(y,z)\in X,} \hfill \\
 { - {\displaystyle\frac{{1}}{{2}}},\,\,\,\,(y,z)\in\gamma,} \hfill \\
 {\,\,\,\,\,\,0,\,\,\,\,\,(y,z)\notin {X\cup \gamma}.} \hfill \\
\end{array}}}  \right.
\end{equation*}
Лемма \ref{L1} доказана.

\begin{lemma} \label{L2}{Если поверхность $\Gamma$ удовлетворяет перечисленным выше ус\-ло\-ви\-ям, то}
     \begin{equation*} \label{eq317}
{\int\int_{\Gamma}{\left|B_\nu^{\alpha}[q_{1} \left(
{\xi,\eta,\zeta;x,y,z} \right)]\right|d\theta d\vartheta}} \leq B,
\end{equation*}
где $B$ -- постоянная.
\end{lemma}

 \textbf{Доказательство.}  Проведем плоскость $x=\delta$ ($\delta>0$ достаточно мало) и части поверхности $\Gamma$, находящиеся ниже и выше этой плоскости обозначим через $H_\delta$ и $L_\delta$, соответственно. Формулу (\ref{eq306}) представим в виде
$$B_\nu^{\alpha}[q_{1} \left( {\xi,\eta,\zeta;x,y,z} \right)]= P\left(\xi,\eta, \zeta; x,y,z\right)+Q\left(\xi,\eta, \zeta; x,y,z\right), $$
где

$P\left(\xi,\eta, \zeta; x,y,z\right)=$
\begin{equation*}
\label{eq1445} =-\frac{1+2\alpha}{2\pi r
r_1^{2+2\alpha}}x\xi^{2\alpha} F\left( {\alpha -
{\frac{{1}}{{2}}}, 1+\alpha ;1 + 2\alpha ;1-\frac{r^2}{r_1^2}}
\right){\cos \left( {\nu,\xi}  \right)},
 \end{equation*}
\begin{equation*}
\label{eq1446} Q\left(\xi,\eta, \zeta; x,y,z\right) =-
\frac{1+2\alpha}{4\pi r r_1^{2\alpha}}F\left( {\alpha -
{\frac{{3}}{{2}}},\alpha ;2\alpha ;1-\frac{r^2}{r_1^2}}
\right)B_\nu^{\alpha}\left[ \ln r^2 \right].
\end{equation*}

Из теории специальных функций известно, что гипергеометрическая
функция Гаусса $F(a,b;c;t)$ при $c-a-b>0$ и $|t|\leq 1$
ограничена. Кроме того, согласно теории потенциала для уравнения
Лапласа имеет место неравенство
 \begin{equation*} \label{eq322}
\left|\int\int_{\Gamma}\mu\frac{\cos(r,\nu)}{r^2}d\theta
d\vartheta\right|<C,
\end{equation*}
где $\mu$ -- ограниченная интегрируемая функция, а $C$ --
некоторая постоянная.

Теперь, легко видеть, что
 \begin{equation} \label{eq1447}
{\int\int_{L_\delta}{\left|P\left({\xi,\eta,\zeta;x,y,z}\right)\right|dL_\delta}}
\leq C_1\,\,\,(\delta>0),
\end{equation}
где $C_1$ не зависит от $(x,y,z).$

В силу  неравенства $r_1^2>r^2$ имеем
\begin{equation}\label{eq1448}
 {\int\int_{H_\delta}{\left|P\left({\xi,\eta,\zeta;x,y,z}\right)\right|dH_\delta}}\leq C_2{\int\int_{H_\delta}\frac{\left|\cos(\nu,\xi)\right|}{r^2}dH_\delta} \leq C_3.
\end{equation}
Далее, проводив аналогичные рассуждения, получим

\begin{equation*}\label{eq1450}
 \int\int_{\Gamma}{\left|Q\right|d\Gamma}\leq
 C_4\int\int_{\Gamma}{\frac{1}{rr_1^{2\alpha}}\left|B_\nu^{\alpha}\left[ \ln {r} \right]\right|d\theta d\vartheta}\leq
\end{equation*}
\begin{equation}\label{eq1450}
 \leq
 C_5{\int\int_{\Gamma}\frac{\left|\cos(\nu,\xi)\right|}{r^2}d\theta d\vartheta} \leq C_6.
\end{equation}
Таким образом, из полученных оценок (\ref{eq1447}),
(\ref{eq1448}), (\ref{eq1450}) следует спра\-вед\-ли\-вость леммы
\ref{L2}.

\begin{lemma} \label{L3}{Если точка $(x,y,z)$ лежит на  $\Gamma$, то }
     \begin{equation}  \label{eq323}
\left|B_\nu^{\alpha}\left[{{{q_{1}
\left({\xi,\eta,\zeta;x,y,z}\right)}}}\right]\right|\leq
\frac{B_1}{r_1^{2\alpha}r}
\end{equation}
где $B_1$ -- постоянная.
\end{lemma}

 \textbf{Доказательство.}
Оценка непосредственно следует из формулы (\ref{eq306}).

Формулы (\ref{eq304}) показывают, что при $\mu_1\equiv1$ потенциал
двойного слоя испытывает разрыв непрерывности, когда точка
$(x,y,z)$ пересекает поверхность $\Gamma$. В случае произвольной
непрерывной плотности $\mu_1(s,t)$  имеет место
\begin{theorem}  \label{T1}
      Потенциал двойного слоя $w^{(1)}(x,y,z)$   имеет пределы при стремлении точки $(x,y,z)$ к точке $(x_0,y_0,z_0)$ поверхности $\Gamma$  извне или изнутри. Если предел значений $w^{(1)}(x,y,z)$ изнутри обозначить через $w_i^{(1)}(s,t)$ , а предел извне  -- через $w_e^{(1)}(s,t)$, то имеют место формулы
\begin{equation}\label{eq324}
 \begin{array}{*{20}c}
 \displaystyle{ w_i^{(1)}(s,t)=-\frac{1}{2}\mu_1(s,t)+\int\int_\Gamma\mu_1(\theta,\vartheta)K_1(s,t;\theta,\vartheta)d\theta d\vartheta,} \hfill \\
 \\
 \displaystyle{ w_e^{(1)}(s,t)=\frac{1}{2}\mu_1(s,t)+\int\int_\Gamma\mu_1(\theta,\vartheta)K_1(s,t;\theta,\vartheta)d\theta d\vartheta,} \hfill \\
 \end{array}
\end{equation}
где
\begin{equation*}\label{eq325}
 K_1(s,t;\theta,\vartheta)=B_\nu^{\alpha}\left[q_1\left(\xi(\theta,\vartheta),\eta(\theta,\vartheta),\zeta(\theta,\vartheta); x(s,t),y(s,t),z(s,t)\right)\right].
\end{equation*}
\end{theorem}

  \textbf{Доказательство}    теоремы  \ref{T1} следует из лемм \ref{L1} и \ref{L2}.

Функция $$w_0^{(1)}(s,t)
=\int\int_\Gamma\mu_1(\theta,\vartheta)K_1(s,t;\theta,\vartheta)d\theta
d\vartheta$$ непрерывна при    $(x,y,z)\in\overline{\Gamma}$, что
следует из хода доказательства теоремы \ref{T1}. Принимая во
внимание формулы (\ref{eq324}) и непрерывность функций
$w_0^{(1)}(s,t)$  и $\mu_1(s,t)$ при
$(s,t)\in[s_1,s_2]\times[t_1,t_2]$, мы можем утверждать, что
потенциал двойного слоя $w^{(1)}(x,y,z)$   есть функция
непрерывная внутри области $D$ вплоть до поверхности $\Gamma$.
Точно также $w^{(1)}(x,y,z)$  непрерывна вне области $D$ вплоть до
поверхности $\Gamma$.

\subsection{Потенциал двойного слоя $w^{\left( {2} \right)}\left( {x},y,z \right)$}\label{s3.2}
Используя второе фун\-да\-мен\-таль\-ное решение (\ref{eq3})
уравнения (\ref{eq1}), определим потенциал двой\-но\-го слоя
формулой
\begin{equation}
\label{eq326} w^{\left( {2} \right)}\left( {x},y,z \right) =
{\int\int_{\Gamma}  {\mu _{2} \left( \theta,\vartheta
\right){{{B_\nu^{\alpha}\left[ q_{2} \left( {\xi,\eta,\zeta;x,y,z}
\right)\right]}}}d\theta d\vartheta}}  ,
\end{equation}
где $\mu _{2} \left( \theta,\vartheta\right)$ -- непрерывная
функция в $\overline{\Gamma}$.

Очевидно, что $w^{\left( {2} \right)}\left( {x},y,z \right)$ есть
регулярное решение уравнения (\ref{eq1}) в любой области, лежащей
в полупространстве $x>0$, не имеющей общих точек ни с поверхностью
$\Gamma$, ни с плоскостью $yOz$.  Потенциал двойного слоя
(\ref{eq326}) оп\-ре\-де\-лен во всех точках полупространства
$x>0$.

\begin{lemma} \label{L5}{Справедлива следующая формула:}
     \begin{equation*} \label{eq328}
w_{1}^{\left( {2} \right)} \left( {x},y,z \right) \equiv
{\int\int_{\Gamma}B_\nu^{\alpha}[q_{2} \left(
{\xi,\eta,\zeta;x,y,z} \right)]d\Gamma} = \left\{
{{\begin{array}{*{20}c}
 { i(x,y,z) - 1,\,\,\,\,x \in D,} \hfill \\
 { i(x,y,z) - {\displaystyle\frac{{1}}{{2}}},\,\,\,\,x \in {\Gamma},} \hfill \\
 { i(x,y,z),\,\,\,\,\,\,\,\,\,\,\,\,\,\,\,\,x \notin {{D \cup \Gamma}},} \hfill \\
\end{array}}}  \right.
\end{equation*}
где
\begin{equation*}\label{eql1333}
i(x,y,z)\equiv{\int\int_{X}\frac{\partial q_{2} \left(
{0,\eta,\zeta;x,y,z} \right)}{\partial \xi}d\eta d\zeta}=
 \end{equation*}
\begin{equation*}\label{eql1334}
=(1-2\alpha) x^{1-2\alpha}{\int\int_{X}\frac{d\eta
d\zeta}{\left[x^2+(y-\eta)^2+(z-\zeta)^2\right]^{3/2-\alpha}}}.
 \end{equation*}
\end{lemma}

Доказательство этой леммы проводится так же, как и доказательство
леммы \ref{L1}.

\begin{lemma} \label{L6}{При любом положении точки $(x,y,z)$ в полупространстве $x>0$ имеет место неравенство }
     \begin{equation*} \label{eq329}
{\int\int_{\Gamma}{\left|B_\nu^{\alpha}[q_{2} \left(
{\xi,\eta,\zeta;x,y,z} \right)]\right|d\Gamma}} \leq B_2,
\end{equation*}
где $B_2$ -- постоянная.
\end{lemma}

\begin{lemma} \label{L7}{Если точка $(x,y,z)$ лежит на  $\Gamma$, то }
     \begin{equation}  \label{eq330}
\left|B_\nu^{\alpha}\left[{{{q_{2}
\left({\xi,\eta,\zeta;x,y,z}\right)}}}\right]\right|\leq
\frac{B_3}{r_1^{2\alpha}r}
\end{equation}
где $B_3$ -- постоянная.
\end{lemma}

Доказательство лемм \ref{L6} и \ref{L7} непосредственно следует из
формулы

$B_\nu^{\alpha}[q_{2} \left( {\xi,\eta,\zeta;x,y,z} \right)]=$
\begin{equation*}
 = -\frac{3-2\alpha}{2\pi}
\frac{x^{1-2\alpha}\xi^{1-2\alpha}}{2r^{3-2\alpha}}F\left(\frac{5}{2}-\alpha,1-\alpha
;2-2\alpha ;1-\frac{r^2}{r_1^2} \right)B_\nu^{\alpha}\left[\ln
r^2\right]
\end{equation*}

\begin{equation*}
\label{eq331}
+\frac{1-2\alpha}{2\pi}\frac{x^{1-2\alpha}}{r^{3-2\alpha}}F\left(\frac{3}{2}-\alpha,1-\alpha
;1-2\alpha ;1-\frac{r^2}{r_1^2} \right)\cos\left( {\nu,\xi}
\right).
 \end{equation*}

\begin{theorem}  \label{T2}
      Для непрерывной плотности $\mu_2(s,t)$  имеют место формулы
\begin{equation}\label{eq332}
 \begin{array}{*{20}c}
 \displaystyle{ w_i^{(2)}(s,t)=-\frac{1}{2}\mu_2(s,t)+\int\int_\Gamma\mu_2(\theta,\vartheta)K_2(s,t;\theta,\vartheta)d\theta d\vartheta,} \hfill \\
 \\
 \displaystyle{ w_e^{(2)}(s,t)=\frac{1}{2}\mu_2(s,t)+\int\int_\Gamma\mu_2(\theta,\vartheta)K_2(s,t;\theta,\vartheta)d\theta d\vartheta,} \hfill \\
 \end{array}
\end{equation}
где
\begin{equation*}\label{eq333}
 K_2(s,t;\theta,\vartheta)=B_\nu^{\alpha}\left[q_2\left(\xi(\theta,\vartheta),\eta(\theta,\vartheta),\zeta(\theta,\vartheta); x(s,t),y(s,t),z(s,t)\right)\right],
 \end{equation*}
точки
$\left(\xi(\theta,\vartheta),\eta(\theta,\vartheta),\zeta(\theta,\vartheta\right)$
и $\left(x(s,t),y(s,t),z(s,t)\right)$ лежат на поверхности
$\Gamma$.
\end{theorem}

Доказательство теоремы \ref{T2} следует из лемм \ref{L5} и
\ref{L6}.

\section{Потенциал простого слоя}  \label{S4}

\subsection{Потенциал простого слоя $v_1( x,y,z)$}\label{s4.1} Пусть поверхность $\Gamma$ удов\-ле\-творяет тем
же условиям, что в \S\, \ref{S3}.

\begin{definition}\label{de2}
Потенциалом простого слоя с плотностью $\rho_1
\left(\theta,\vartheta \right)$ назовем функцию
\begin{equation}
\label{eq401} v_1\left( x,y,z \right) = {\int\int_{\Gamma}  {
\rho_1 \left(\theta,\vartheta \right){q_{1} \left(
{\xi,\eta,\zeta;x,y,z} \right)}d\theta d\vartheta}}  ,
\end{equation}
где $q_{1} \left( {\xi,\eta,\zeta;x,y,z} \right)$  --
фундаментальное решение уравнения (\ref{eq1}). \end{definition}
Будем предполагать, что $\rho_1 \left( \theta,\vartheta\right)$ --
непрерывная функция на $\overline{\Gamma}$. Потенциал простого
слоя (\ref{eq401}) определен во всем полупространстве $x>0$ и
остается непрерывным при переходе через поверхность $\Gamma$.
Очевидно, что потенциал простого слоя $v_1\left( x,y,z\right)$
есть регулярное решение уравнения (\ref{eq1}) в любой области,
лежащей в полупространстве $x>0$, не имеющей общих точек ни с
поверхностью $\Gamma$, ни с плоскостью $yOz$. Нетрудно видеть, что
при стремлении точки $(x,y,z)$  к бесконечности потенциал простого
слоя  $v_1\left( x,y,z\right)$  стремится к нулю. Действительно,
пусть точка $(x,y,z)$ находится на полусфере $C_R$:
$x^2+y^2+z^2=R^2\,(x>0)$, тогда в силу (\ref{eq2}) имеем
\begin{equation}
\label{eq402} \left|v_1\left( x,y,z \right)\right| \leq
{\int\int_{\Gamma}  {| \rho_1 \left(
\theta,\vartheta\right)||{q_{1} \left( {\xi,\eta,\zeta;x,y,z}
\right)}|d\Gamma}} \leq MR^{-1-2\alpha} \,\,\left(R\geq
R_0\right),
\end{equation}
где $M$ -- постоянная.

\subsection{Конормальная производная потенциала простого слоя}\label{s4.2} Возь\-мем на поверхности $\Gamma$
произвольную точку $N\left(x(s,t),y(s,t),z(s,t)\right)$ и проведем
в этой точке нормаль. Рассмотрим на этой нормали какую-нибудь
точку $M\left(x,y,z\right)$, не лежащую на поверхности $\Gamma$, и
составим конормальную производную от потенциала простого слоя
(\ref{eq401}):
\begin{equation}
\label{eq403} B_n^{\alpha}\left[{v_1\left( x,y,z \right)}\right] =
{\int\int_{\Gamma}  {\rho_1 \left(
\theta,\vartheta\right)B_n^{\alpha}\left[{q_1\left(\xi,\eta,\zeta;
x,y,z \right)}\right]d\theta d\vartheta}},
\end{equation}
где
\begin{equation*}
B_n^{\alpha}[\,\,\,]=x^{2\alpha}\left(cos(n,x)\cdot\frac{\partial}{\partial
x}+ cos(n,y)\cdot\frac{\partial}{\partial y}+
cos(n,z)\cdot\frac{\partial}{\partial z}\right).
\end{equation*}

Интеграл (\ref{eq403}) существуют и в том случае, когда точка
$M(x,y,z)$  совпадает с точкой  $N\left(x_0,y_0,z_0\right)$,
упомянутой выше.

Обозначим через  $B_n^{\alpha}\left[{v_1\left( x,y,z
\right)}\right]_i$  и $B_n^{\alpha}\left[{v_1\left( x,y,z
\right)}\right]_e$  соответственно пре\-дель\-ные значения
нормальной производной при стремлении точки $M(x,y,z)$ к точке
$N\in \Gamma$ изнутри и извне поверхности $\Gamma$.
\begin{theorem} \label{T3}
Для непрерывной плотности $\rho_{1}(\xi,\eta,\zeta)$ имеют место
сле\-ду\-ю\-щие формулы:
\begin{equation}\label{eq404}
\begin{array}{*{20}c}
 \displaystyle{ B_n^{\alpha}\left[{v_1\left( x,y,z \right)}\right]_i=\frac{1}{2}\rho_1(s,t)+\int\int_\Gamma\rho_1(\theta,\vartheta)K_1(\theta,\vartheta; s,t)d\theta d\vartheta,} \hfill \\
 \\
\displaystyle{B_n^{\alpha}\left[{v_1\left( x,y,z \right)}\right]_e=-\frac{1}{2}\rho_1(s,t)+\int\int_\Gamma\rho_1(\theta,\vartheta)K_1(\theta,\vartheta; s,t)d\theta d\vartheta,} \hfill \\
 \end{array}
\end{equation}
где
\begin{equation*}\label{eq405}
K_1(\theta,\vartheta; s,t)=B_n^{\alpha}\left[{q_1\left(
\xi(\theta,\vartheta),\eta(\theta,\vartheta),\zeta(\theta,\vartheta);
x(s,t),y(s,t),z(s,t) \right)}\right].
\end{equation*}
\end{theorem}

Из этих формул непосредственно следует величина скачка нормальной
производной потенциала простого слоя:
\begin{equation}\label{eq406}
  B_n^{\alpha}\left[{v_1\left( x,y,z \right)}\right]_i-B_n^{\alpha}\left[{v_1\left( x,y,z \right)}\right]_e=\rho_1(x,y,z)
\end{equation}

Совершенно так же, как и в неравенстве (\ref{eq402}), можно
показать, что имеет место следующая оценка:
\begin{equation}
\label{eq407} \left|B_n^{\alpha}\left[{v_1\left( x,y,z
\right)}\right]\right|
 \leq MR^{-2-2\alpha} \,\,\left(R\geq R_0\right),
\end{equation}
где $M$  -- постоянная.

\subsection{Применимость формулы Грина для потенциалов}\label{s4.3} Покажем, что для потенциала простого слоя
(\ref{eq401}) применима формула Грина (\ref{eq10}). Рассмотрим
область $D_{\varepsilon,\delta}$, лежащую внутри $D$ и
ограниченную поверхностью $\Gamma_\varepsilon$ параллельной
поверхности $\Gamma$ и замкнутой областью $X_\delta$, т.е.
сечением области $D_{\varepsilon,\delta}$ плоскостью $x=\delta$.
Применим формулу Грина (\ref{eq10}) к потенциалу простого слоя
$v_1(x,y,z)$, выбирая за область интегрирования область
$D_{\varepsilon,\delta}$.  Тогда получим
\begin{equation*}
\label{eq1001} {\int\int\int_{D_{\varepsilon,\delta}} {x^{2\alpha}
\left[{{\left( {{\frac{{\partial v_1}}{{\partial x}} }}
\right)}}^{2}+{{\left( {{\frac{{\partial v_1}}{{\partial y}} }}
\right)}}^{2}+{{\left( {{\frac{{\partial v_1}}{{\partial z}} }}
\right)}}^{2}\right]dxdydz}}
\end{equation*}
\begin{equation}
\label{eq408}  = {\int\int_{\Gamma_\varepsilon} {v_1
B_n^\alpha\left[{v_1}\right]d\Gamma_\varepsilon}}+{\int\int_{X_\delta}
{\delta^{2\alpha} v_1\left(\delta,y,z\right){\frac{{\partial
v_1\left(\delta,y,z\right)}}{{\partial x}}}dydz}} .
\end{equation}

Конормальная производная  $B_n^{\alpha}\left[{v_1}(x,y,z)\right]$
есть непрерывная функция вплоть до поверхности $\Gamma$ и
$${\mathop {\lim} \limits_{x \to
0}}x^{2\alpha}{\displaystyle\frac{{\partial
v_1\left(x,y,z\right)}}{{\partial x}}}=0, \,\, (y,z)\in X.$$
Следовательно,
\begin{equation}
\label{eq409}   {\mathop {\lim} \limits_{\delta \to
0}}{\int\int_{X_\delta} {\delta^{2\alpha}
v_1\left(\delta,y,z\right){\frac{{\partial
v_1\left(\delta,y,z\right)}}{{\partial x}}}dydz}} = 0.
\end{equation}

Переходя в формуле  (\ref{eq408}) к пределу  при  $\varepsilon\to
0$, $\delta\to 0$  и принимая во внимание  (\ref{eq409}), получим
\begin{equation*}
\label{eq412} {\int\int\int_{D}  {x^{2\alpha} \left[{{\left(
{{\frac{{\partial v_1}}{{\partial x}} }} \right)}}^{2}+{{\left(
{{\frac{{\partial v_1}}{{\partial y}} }} \right)}}^{2}+{{\left(
{{\frac{{\partial v_1}}{{\partial z}} }}
\right)}}^{2}\right]dxdydz}}=
\end{equation*}
\begin{equation}
\label{eq412} ={\int\int_{\Gamma} v_1B_n^{\alpha}\left[{v_1\left(
x,y,z \right)}\right]_id\Gamma}.
\end{equation}

Плоскую область, ограниченную кривой $\gamma$ и окружностью
$y^2+z^2=R^2$, лежащей в плоскости $yOz$, обозначим через
$X_{0R}$. Применим теперь формулу (\ref{eq412}) к области $D_R'$,
ограниченной поверхностью $\Gamma$, областью $X_{0R}$ и полусферой
$C_R$, содержащей область $D$. Переходя затем к пределу при $R \to
\infty$ и учитывая  (\ref{eq402}) и (\ref{eq407}), получим
\begin{equation*}
 {\int\int\int_{D'}  {x^{2\alpha} \left[{{\left(
{{\frac{{\partial v_1}}{{\partial x}} }} \right)}}^{2}+{{\left(
{{\frac{{\partial v_1}}{{\partial y}} }} \right)}}^{2}+{{\left(
{{\frac{{\partial v_1}}{{\partial z}} }}
\right)}}^{2}\right]dxdydz}}=
\end{equation*}
\begin{equation}
\label{eq413} =-{\int\int_{\Gamma}
{v_1B_n^{\alpha}\left[{v_1\left( x,y,z \right)}\right]_ed\Gamma}}.
\end{equation}
Здесь и далее $D'=R_3^+\setminus\bar{D}$ -- неограниченная область
при $x>0$.

\subsection{Потенциал простого слоя $v_2(x,y,z)$} \label{s4.4} Используя второе фун\-да\-мен\-таль\-ное решение
(\ref{eq3}) уравнения (\ref{eq1}), определим потенциал
про\-сто\-го слоя формулой
\begin{equation}
\label{eq414} v_2(x,y,z) = {\int\int_{\Gamma}  {\rho _{2} \left(
\theta,\vartheta\right){{{ q_{2} \left( {\xi,\eta,\zeta;x,y,z}
\right)}}}d\theta d\vartheta}}  ,
\end{equation}
где $\rho_{2} \left(\xi,\eta,\zeta \right)$ -- не\-пре\-рыв\-ная
функция в $\overline{\Gamma}$. Потенциал простого слоя
(\ref{eq414}) не\-пре\-ры\-вен во всем полупространстве $x>0$.
Очевидно, что $v_2( {x},y,z )$ есть регулярное решение уравнения
(\ref{eq1}) в любой области, лежащей в по\-лу\-про\-ст\-ранстве
$x>0$, не имеющей общих точек ни с поверхностью $\Gamma$, ни с
плоскостью $yOz$.

Конормальная производная от потенциала простого  слоя
(\ref{eq414}) равна
\begin{equation*}
\label{eq415} B_n^\alpha \left[v_2( {x},y,z)\right] =
{\int\int_{\Gamma}  \rho _{2} \left(
\theta,\vartheta\right)B_n^\alpha\left[{{{ q_{2} \left(
{\xi,\eta,\zeta;x,y,z} \right)}}}\right]d\theta
d\vartheta},\,\,(x,y,z)\notin\Gamma.
\end{equation*}

\begin{theorem} \label{T4}
Для непрерывной плотности $\rho_{2}(\theta,\vartheta)$ имеют место
сле\-ду\-ю\-щие формулы:
\begin{equation}\label{eq416}
\begin{array}{*{20}c}
 \displaystyle{ B_n^{\alpha}\left[{v_2\left( x,y,z \right)}\right]_i=\frac{1}{2}\rho_2(s,t)+\int\int_\Gamma\rho_2(\theta,\vartheta)K_2(\theta,\vartheta; s,t)d\theta d\vartheta,} \hfill \\
 \\
\displaystyle{B_n^{\alpha}\left[{v_2\left( x,y,z \right)}\right]_e=-\frac{1}{2}\rho_2(s,t)+\int\int_\Gamma\rho_2(\theta,\vartheta)K_2(\theta,\vartheta; s,t)d\theta d\vartheta,} \hfill \\
 \end{array}
\end{equation}
где
\begin{equation*}\label{eq417}
K_2(\theta,\vartheta; s,t)=B_n^{\alpha}\left[{q_2\left(
\xi(\theta,\vartheta),\eta(\theta,\vartheta),\zeta(\theta,\vartheta);
x(s,t),y(s,t),z(s,t) \right)}\right].
\end{equation*}
\end{theorem}
Из этих формул непосредственно следует величина скачка нормальной
про\-из\-вод\-ной потенциала простого слоя:
\begin{equation*}\label{eq418}
  B_n^{\alpha}\left[{v_2\left( x,y,z \right)}\right]_i-B_n^{\alpha}\left[{v_2\left( x,y,z \right)}\right]_e=\rho_2(x,y,z)
\end{equation*}

Для дальнейшего полезно заметить, что при стремлении точки
$(x,y,z)$ к бесконечности имеют место оценки
\begin{equation*}
\label{eq419} \left|v_2(x,y,z)\right|\leq \frac{M}{R}, \,\,\,
\left|B_n^{\alpha}\left[{v_2\left( x,y,z \right)}\right]\right|
 \leq MR^{-2-2\alpha} \,\,\left(R\geq R_0\right),
\end{equation*}
где точка $(x,y,z)$ находится на полусфере $C_R:$
$x^2+y^2+z^2=R^2,\,$ $M$  -- постоянная.

Совершенно так же, как это было для потенциала простого слоя
$v_1(x,y,z)$, можно показать, что для потенциала простого слоя
(\ref{eq414}) применимы формулы Грина  (\ref{eq412}) и
(\ref{eq413}).

\section{Интегральные уравнения для плотностей} \label{S5}

Формулы (\ref{eq324})  и (\ref{eq404}), а также (\ref{eq332}) и
(\ref{eq416}) могут быть написаны как интегральные уравнения для
плотностей:
\begin{equation}\label{eq501}
\mu_j(s,t)-\lambda\int\int_\Gamma
K_j(s,t;\theta,\vartheta)\mu_j(\theta,\vartheta)d\theta
d\vartheta=f_j(s,t)\,\,\, \left(j=1.2\right),
\end{equation}

\begin{equation}\label{eq502}
\rho_j(s,t)-\lambda\int\int_\Gamma
K_j(\theta,\vartheta;s,t)\rho_j(\theta,\vartheta)d\theta
d\vartheta=g_j(s,t)\,\,\, \left(j=1.2\right),
\end{equation}
где
\begin{equation*}
\lambda=2,\,\,\,\,f_j(s,t)=-2w_i^{(j)}(s,t),\,\,
g_j(s,t)=-2B_n^\alpha\left[v_j\right]_e,
\end{equation*}

\begin{equation*}
\lambda=-2,\,\,\,\,f_j(s,t)=2w_e^{(j)}(s,t),\,\,
g_j(s,t)=2B_n^\alpha\left[v_j\right]_i.
\end{equation*}

Уравнения (\ref{eq501}) и (\ref{eq502}) сопряженные и в силу леммы
\ref{L3} к ним применима теория Фредгольма.

Покажем, что $\lambda=2$ не является собственным значением ядра
$K_1(s,t;\theta,\vartheta)$ . Это утверждение эквивалентно тому,
что однородное интегральное уравнение
\begin{equation}\label{eq503}
\rho(s,t)-2\int\int_\Gamma
K_1(\theta,\vartheta;s,t)\rho(\theta,\vartheta)d\theta
d\vartheta=0,
\end{equation}
не имеет нетривиальных решений. Пусть  $\rho_0(\theta,\vartheta)$
-- непрерывное нетривиальное решение уравнения (\ref{eq503}).
Потенциал простого слоя с плотностью $\rho_0(\theta,\vartheta)$
даст нам функцию $v_0(x,y,z)$, которая является решением уравнения
(\ref{eq1}) в областях $D$ и $D'$ и у которой предельные значения
конормальной производной $B_n^\alpha\left[v_0\right]_e$  равны
нулю в силу уравнения (\ref{eq503}). К потенциалу простого слоя
$v_0(x,y,z)$  применима формула (\ref{eq413}), из которой следует,
что $v_0(x,y,z)=const$ в области $D'$. На бесконечности потенциал
простого слоя равен нулю и, следовательно, $v_0(x,y,z)\equiv 0$  в
$D'$, а также и на поверхности $\Gamma$. Применяя теперь формулу
(\ref{eq412}), мы получим, что $v_0(x,y,z)\equiv 0$ и внутри
области $D$. Но тогда $B_n^\alpha \left[v_0\right]_i=0$, и на
основании формулы (\ref{eq406}) получим
$\rho_0(\theta,\vartheta)\equiv 0$. Таким образом, однородное
уравнение (\ref{eq503}) имеет только тривиальное решение;
следовательно, $\lambda=2$ не есть собственное значение ядра
$K_1(s,t;\theta,\vartheta)$.

Однородное уравнение
\begin{equation*}\label{eq504}
\mu(s,t)-\lambda\int\int_\Gamma
K_1(s,t;\theta,\vartheta)\mu(\theta,\vartheta)d\theta
d\vartheta=0,
\end{equation*}
при $\lambda=-2$ имеет в силу (\ref{eq304}) решение, равное
произвольной постоянной, т.е. $\lambda=-2$  есть собственное
значение ядра $K_1(s,t;\theta,\vartheta)$. Совершенно так же можно
показать, что $\lambda=2$ и $\lambda=-2$  не являются собственными
значениями ядра $K_2(s,t;\theta,\vartheta)$.

\section{Решение краевых задач с помощью потенциалов} \label{S6}

\subsection{Постановка краевых задач Дирихле и Хольмгрена}\label{s6.0} Пусть $D$ -- область, ограниченная
односвязной открытой областью $X$ плоскости $yOz$ и поверхностью
$\Gamma$,  лежащей в полупространстве $x>0$. Общую границу плоской
области $X$ и поверхности $\Gamma$ обозначим через $\gamma$. Будем
предполагать, что поверхность $\Gamma$ удовлетворяет условиям 1) и
2) \S \ref{S3}. Рассмотрим две краевые задачи для уравнения
(\ref{eq1}).

\textbf{Задача Дирихле}. Найти в области $D$ регулярное решение
уравнения (\ref{eq1}), непрерывное в замкнутой области
$\overline{D}$ и удовлетворяющее краевому условию
\begin{equation}\label{eq601}
  \left.u\right|_\Gamma=\varphi(s,t),\,\,(s,t)\in \overline{\Phi}; \,\,\, u(0,y,z)=\tau_1(y,z),\,\, (y,z)\in\overline{X},
\end{equation}
где $\varphi(s,t)$ и $\tau_1(y,z)$  -- заданные непрерывные
функции, причем  $ \left.\varphi(s,t)\right|_\gamma=
\left.\tau_1(y,z)\right|_\gamma$.
\\

\textbf{Задача Хольмгрена}. Найти в области $D$ регулярное решение
уравнения (\ref{eq1}), непрерывное в замкнутой области
$\overline{D}$ и удовлетворяющее краевым условиям
\begin{equation}\label{eq602}
  \left.u\right|_\Gamma=\varphi(s,t),\,\,(s,t)\in \overline{\Phi}; \,\,\, \lim\limits_{x \to 0}\left(x^{2\alpha}\frac{\partial u(x,y,z)}{\partial x}\right)=\nu_1(y,z),\,\, (y,z)\in {X},
\end{equation}
где $\nu_1(y,z)$ -- непрерывная функция в $X$, причем в случае
когда кривая $\gamma$ есть окружность, при стремлении точек
$(y,z)\in X$ к кривой $\gamma$ функция $\nu_1(y,z)$ может
обращаться в бесконечность порядка меньше $1-2\alpha$.
\\

\subsection{Единственность решения краевых задач Дирихле и Хольм\-гре\-на}\label{s6.01}

Рассмотрим область $D_{\varepsilon,\delta}$, лежащую внутри $D$ и
ограниченную поверхностью $\Gamma_\varepsilon$ параллельной
поверхности $\Gamma$ и замкнутой областью $X_\delta$, т.е.
сечением области $D_{\varepsilon,\delta}$ плоскостью $x=\delta$.
 Интегрируя обе части тождества
(\ref{eq7}) по области $D_{\varepsilon,\delta}$, и пользуясь
формулой Гаусса-Остроградского, по\-лу\-чим

\begin{equation}
\label{eqq8} \int\int\int_{D_{\varepsilon,\delta}}  x^{2\alpha}
\left[ uE(v)-vE(u) \right]dxdydz =
\int\int_{S_{\varepsilon,\delta}}
{{\left({uB_n^{\alpha}[v]-vB_n^{\alpha}[u]}\right)}}dS_{\varepsilon,\delta},
\end{equation}
где $S_{\varepsilon,\delta} $ граница области
$D_{\varepsilon,\delta} $.

Нетрудно убедиться в справедливости следующего равенства:

\begin{equation*}
{\int\int\int_{D_{\varepsilon,\delta} }   {x^{2\alpha}
uE(u)dxdydz}} = {\int\int\int_{D_{\varepsilon,\delta} }
{x^{2\alpha}\left[u_x^2+u_y^2+u_z^2\right]}dxxydz}-
\end{equation*}
\begin{equation*}-
 \int\int\int_{D_{\varepsilon,\delta}}\left[\frac{\partial}{\partial x}\left(x^{2\alpha}uu_x\right)
 +x^{2\alpha}\frac{\partial}{\partial y}\left(uu_y\right)+x^{2\alpha}\frac{\partial}{\partial z}\left(uu_z\right)\right] dxdydz .
\end{equation*}

Применяя формулу Гаусса-Остроградского к этому тождеству, после
перехода к пределу при  $\delta \to 0$ и $\varepsilon \to 0$,
имеем
\begin{equation*}
\label{eqq41} {\int\int\int_{D}
{x^{2\alpha}\left[u_x^2+u_y^2+u_z^2\right]dxdydz}}  =
\end{equation*}
\begin{equation}
\label{eqq41}   = {\int\int_{X} {\tau_1(y,z)\nu_1 (y,z)dydz}} -
\int\int_{\Gamma}\varphi(x,y,z) B_{nx}^\alpha[u] dxdydz.
\end{equation}

Если теперь рассмотрим однородную задачу Дирихле (Хольмгрена), то
из (\ref{eqq41}) получим

\begin{equation*}
{\int\int\int_{D}  {x^{2\alpha}
\left[u_x^2+u_y^2+u_z^2\right]dxdydz}} = 0.
\end{equation*}
Отсюда следует, что  $u(x,y,z) = 0 $ в  $\overline {D} $.

Справедлива следующая

\begin{theorem} \label{Tedinst} Если задача Дирихле (Хольмгрена) для уравнения (\ref{eq1}) имеет регулярное решение,  то оно единственно.
\end{theorem}

\subsection{Сведение краевых задач к случаю однородных граничных условий на плоскости вырождения}\label{s6.1}

Покажем сначала, что можно ограничиться случаем, когда
$\tau_1(y,z)\equiv 0$  и соответственно $\nu_1(y,z)\equiv 0$.

Пусть $X=\left\{(y,z): p< y < q,\,\, h(y)<z<k(y), \,\right\}$,\,
где  $p$ и $q$ -- некоторые конечные действительные числа, а
$z=h(y)$ и $z=k(y)$  -- непрерывные функции на отрезке $[p,\,q]$.

В случае задачи Хольмгрена, используя фундаментальное решение
(\ref{eq2}) уравнения (\ref{eq1}), возьмем решение в виде
\begin{equation}\label{eq603}
v_1(x,y,z)=-\frac{1}{2\pi}\int\int_{X}\nu_1(\eta,\zeta)\left[x^2+(\eta-y)^2+(\zeta-z)^2\right]^{-\frac{1}{2}-\alpha}d\eta
d\zeta.
\end{equation}

Дифференцируя по  $x$  интеграл (\ref{eq603}), получим
\begin{equation*}\label{eq605}
\frac{\partial v_1}{\partial x}
=\frac{1+2\alpha}{2\pi}x\int\int_{X}\nu_1(\eta,\zeta)\left[x^2+(\eta-y)^2+(\zeta-z)^2\right]^{-\frac{3}{2}-\alpha}d\eta
d\zeta.
\end{equation*}

Нетрудно видеть, что
\begin{equation}\label{eq606}
  {\lim\limits_{{{\begin{smallmatrix}
 {x \to 0 } \\
 {y \to y_0} \\
 {z \to z_0}
\end{smallmatrix}}}}}\frac{1+2\alpha}{2\pi}x^{1+2\alpha}\int\limits_{p}^{q}d\eta\int\limits_{h(\eta)}^{k(\eta)}\left[x^2+(\eta-y)^2+(\zeta-z)^2\right]^{-\frac{3}{2}-\alpha}
d\zeta =1.
\end{equation}
Действительно, полагая $\eta=y+xs$ и $\zeta=z+xt$, будем иметь

\begin{equation*}\frac{1+2\alpha}{2\pi}x^{1+2\alpha}\int\limits_{p}^{q}d\eta\int\limits_{h(\eta)}^{k(\eta)}\nu_1(\eta,\zeta)\left[x^2+(\eta-y)^2+(\zeta-z)^2\right]^{-\frac{3}{2}-\alpha}
d\zeta =
\end{equation*}
\begin{equation*}=\frac{1+2\alpha}{2\pi}\int\limits_{\frac{p-y}{x}}^{\frac{q-y}{x}}\,\,\,\int\limits_{\frac{h(y+xs)-z}{x}}^{\frac{k(y+xs)-z}{x}}\nu_1(y+xs,z+xt)\left[1+s^2+t^2\right]^{-\frac{3}{2}-\alpha}
ds dt.
\end{equation*}
Отсюда, учитывая формулу (\ref{eq316}) при  $x\to 0$, $y \to y_0$,
$z \to z_0$  $\left((y_0,\,z_0)\in X\right)$, получим
(\ref{eq606}).

Далее, следуя \cite[гл.2, \S 6]{Sm} и используя равенство
(\ref{eq606}) можно доказать, что производная
$\displaystyle\frac{\partial v_1}{\partial x}$ с весом
$x^{2\alpha}$ в области $X$ принимает значение $\nu_1(y,z)$, т.е.

\begin{equation*}\label{eq604}
  {\lim\limits_{{{\begin{smallmatrix}
 {x \to 0 } \\
 {y \to y_0} \\
 {z \to z_0}
\end{smallmatrix}}}}} \left(x^{2\alpha}\frac{\partial v_1(x,y,z)}{\partial x}\right)= \nu_1\left(y_0,z_0\right), \,\,\,\left(y_0,z_0\right)\in X.
\end{equation*}

Нетрудно видеть, что значения функции $v_1(x,y,z)$ непрерывны и
ограничены на поверхности $\Gamma$. Теперь  решение первоначальной
задачи    Хольмгрена может быть представлено в виде

\begin{equation*}\label{eq609}
  u(x,y,z)=v_1(x,y,z)+w_1(x,y,z),
\end{equation*}
где $w_1(x,y,z)$ есть решение уравнения (\ref{eq1}) в области $D$,
удовлетворяющее краевым условиям
\begin{equation}\label{eq610}
  \left. w_1\right|_\Gamma=\varphi(s,t)-\left. v_1\right|_\Gamma=\varphi_1(s,t), \,\,\,\left.\left(x^{2\alpha}\frac{\partial w_1}{\partial x}\right) \right|_{x=0}=0,\,\,(y,z)\in X.
\end{equation}

Итак, доказано, что в случае задачи Хольмгрена можно ограничится
случаем $\nu_1(y,z)\equiv 0$.

В случае задачи Дирихле выберем  радиус $R$  такой, чтобы область
$D$ целиком лежала в бесконечном цилиндре $y^2+z^2<R^2$. Функцию
$\tau_1(y,z)$ продолжим на весь цилиндр $y^2+z^2\leq R^2$ так,
чтобы продолженная функция  $\tau_1(y,z)$ принадлежала классу
непрерывных функций. Рассмотрим функцию

$v_2(x,y,z)=$
\begin{equation*}\label{eq611}
=\frac{1}{2\pi}x^{1-2\alpha}\int\limits_{-R}^R\int\limits_{-\sqrt{R^2-\zeta^2}}^{\sqrt{R^2-\zeta^2}}\tau_1(\eta,\zeta)\left[x^2+(\eta-y)^2+(\zeta-z)^2\right]^{\alpha-\frac{3}{2}}d\eta
d\zeta.
\end{equation*}
Нетрудно проверить, что функция $v_2(x,y,z)$ является регулярным
решением уравнения  (\ref{eq1}) в полупространстве $x>0$  и
принимает в области $X$ плоскости  $yOz$ значения $\tau_1(y,z)$.
Последнее утверждение доказывается так же, как в случае задачи
Хольмгрена.

Решение первоначальной задачи Дирихле может быть представлено в
виде
\begin{equation*}\label{eq612}
  u(x,y,z)=v_2(x,y,z)+w_2(x,y,z),
\end{equation*}
где $w_2(x,y,z)$ есть решение уравнения (\ref{eq1}) в области $D$,
удовлетворяющее краевым условиям
\begin{equation}\label{eq613}
  \left. w_2\right|_\Gamma=\varphi(s,t)-\left. v_2\right|_\Gamma=\varphi_2(s,t),\,\,\, w_2(0,y,z)=0,\,\,\,(y,z)\in \overline{X}.
\end{equation}

Таким образом, в случае задачи Дирихле можно также ограничиться
случаем $\tau_1(y,z)\equiv 0$.

\subsection{Сведения краевых задач к интегральным уравнениям}\label{s6.2}
Решение задачи Хольмгрена будем искать в виде потенциала двойного
слоя с неизвестной плотностью $\mu_1(\theta,\vartheta)$:
\begin{equation*}\label{eq614}
  w_1(x,y,z)=\int\int_\Gamma\mu_1(\theta,\vartheta)B_\nu^\alpha\left[q_1(\xi,\eta,\zeta;x,y,z)\right]d\theta d\vartheta.
\end{equation*}
Потенциал двойного слоя $\,w_1(x,y,z)\,$ удовлетворяет уравнению
(\ref{eq1}) внутри области $\,D$  и
$\left.\left(x^{2\alpha}\displaystyle\frac{\partial w_1}{\partial
x}\right) \right|_{x=0}=0$ при $(y,z)\in X$. Для выполнения
краевого условия (\ref{eq610}) на $\Gamma$ нужно, чтобы предельные
значения $w_1(x,y,z)$ изнутри равнялись $\varphi_1(s,t)$:
$$
w_{1i}(x,y,z)=\varphi_1(s,t), \,\,\,(x,y,z) \in \Gamma.
$$
Пользуясь первой из формул (\ref{eq324}), получим для плотности
$\mu_1(\theta,\vartheta)$ интегральное уравнение
\begin{equation}\label{eq615}
\mu_1(s,t)-2\int\int_\Gamma
K_1(s,t;\theta,\vartheta)\mu_1(\theta,\vartheta)d\theta
d\vartheta=-2\varphi_1(s,t),
\end{equation}
где
\begin{equation*}\label{eq616}
K_1(s,t;\theta,\vartheta)=B_\nu^{\alpha}\left[{q_1\left(
\xi(\theta,\vartheta),\eta(\theta,\vartheta),\zeta(\theta,\vartheta);
x(s,t),y(s,t),z(s,t) \right)}\right].
\end{equation*}
Как это следует из оценки (\ref{eq323}), ядро
$K_1(s,t;\theta,\vartheta)$ имеет слабую особенность.

Аналогично решение задачи Дирихле будем искать виде потенциала
двойного слоя с неизвестной плотностью $\mu_2(\theta,\vartheta)$:
\begin{equation*}\label{eq614}
  w_2(x,y,z)=\int\int_\Gamma\mu_2(\theta,\vartheta)B_\nu^\alpha\left[q_2(\xi,\eta,\zeta;x,y,z)\right]d\theta d\vartheta.
\end{equation*}

Потенциал двойного слоя $w_2(x,y,z)$ также удовлетворяет уравнению
(\ref{eq1}) внутри области $D$, и $w_2(0,y,z)=0$  при $(y,z)\in
\overline{X}$. Для выполнения краевого условия (\ref{eq613}) на
$\Gamma$ нужно, чтобы предельные значения $w_2(x,y,z)$ изнутри
равнялись   $\varphi_2(s,t)$:
$$
w_{2i}(x,y,z)=\varphi_2(s,t), \,\,\,(x,y,z)\in \Gamma.
$$
Пользуясь первой из формул  (\ref{eq332}), мы приходим к
интегральному уравнению
\begin{equation}\label{eq618}
 \mu_2(s,t)-2\int\int_\Gamma K_2(s,t;\theta,\vartheta)\mu_2(\theta,\vartheta)d\theta d\vartheta=-2\varphi_2(s,t),
\end{equation}
где
\begin{equation*}\label{eq619}
 K_2(s,t;\theta,\vartheta)=B_\nu^{\alpha}\left[q_2\left(\xi(\theta,\vartheta),\eta(\theta,\vartheta),\zeta(\theta,\vartheta); x(s,t),y(s,t),z(s,t)\right)\right].
\end{equation*}
Как это следует из оценки (\ref{eq330}), ядро
$K_2(s,t;\theta,\vartheta)$ имеет слабую особенность.

Отметим, что если какое-либо из уравнений (\ref{eq615}),
(\ref{eq618}) разрешимо, то его решение непрерывно. Это следует из
непрерывности свободного члена и из вида ядер
$K_i(s,t;\theta,\vartheta)$ $(i=1,2)$.

В \S\ref{S5} было показано, что $\lambda=2$  не является
собственным значением ядер $K_1(s,t;\theta,\vartheta)$  и
$K_2(s,t;\theta,\vartheta)$  и, следовательно, уравнения
(\ref{eq615}) и (\ref{eq618}) при любом свободном члене имеют
единственное решение. Таким образом, если поверхность $\Gamma$
удовлетворяет условиям 1) и 2) \S \ref{S3} и заданные значения
решения на поверхности $\Gamma$ непрерывны, то задачи Дирихле и
Хольмгрена для уравнения (\ref{eq1}) имеют единственное решение, и
это решение можно представить в виде потенциала двойного слоя.

\section{Функция Грина оператора E(u)} \label{S7}

\subsection{Функция Грина задачи Хольмгрена}\label{s7.1}
\begin{definition}\label{de3}
Функцией Грина задачи Хольмгрена для уравнения (\ref{eq1})
называется функция $G_1(x,y,z;x_0,y_0,z_0)$, удовлетворяющая
следующим условиям:

1) внутри области $D$, кроме точки $(x_0,y_0,z_0)$, эта функция
есть регулярное решение уравнения (\ref{eq1});

2) она удовлетворяет граничным условиям
\begin{equation}\label{eq701}
  \left. G_1(x,y,z;x_0,y_0,z_0)\right|_\Gamma=0,\,\,\,\left.\left(x^{2\alpha}\frac{\partial G_1}{\partial x}\right)\right|_{x=0}=0;
\end{equation}

3) она может быть представлена в виде
\begin{equation}\label{eq702}
G_1(x,y,z;x_0,y_0,z_0)=q_1(x,y,z;x_0,y_0,z_0)+v_1(x,y,z;x_0,y_0,z_0)
\end{equation}
где $q_{1} \left( {x,y,z;\xi,\eta,\zeta}\right)$  --
фундаментальное решение уравнения (\ref{eq1}), определенное
формулой (\ref{eq2}), а $v_1(x,y,z;x_0,y_0,z_0)$  -- регулярное
решение уравнения (\ref{eq1}) везде внутри $D$. \end{definition}

Построение функции Грина сводится к нахождению ее регулярной части
$v_1(x,y,z;x_0,y_0,z_0)$, которая в силу (\ref{eq5}),
(\ref{eq701}) и (\ref{eq702}) должна удовлетворять граничным
условиям
\begin{equation}\label{eq703}
  \left. v_1(x,y,z;x_0,y_0,z_0)\right|_\Gamma=\left.-q_1(x,y,z;x_0,y_0,z_0)\right|_\Gamma,
\end{equation}
\begin{equation*}\label{eq704}
  \left.\left(x^{2\alpha}\frac{\partial v_1(x,y,z;x_0,y_0,z_0)}{\partial x}\right)\right|_{x=0}=0.
\end{equation*}

Функцию $v_1(x,y,z;x_0,y_0,z_0)$  будем искать в виде потенциала
двойного слоя:
\begin{equation}
\label{eq705} v_1\left( x,y,z;x_0,y_0,z_0 \right) =
{\int\int_{\Gamma}  {\mu _{1} \left( \theta,\vartheta;x_0,y_0,z_0
\right)B_\nu^{\alpha}[q_{1} \left( {\xi,\eta,\zeta;x,y,z}
\right)]d\theta d\vartheta}}.
\end{equation}
Принимая во внимание первое из равенств (\ref{eq324}) и граничное
условие (\ref{eq703}), получим интегральное уравнение для
плотности $\mu_1\left(s,t;x_0,y_0,z_0 \right)$
\begin{equation*}
 \mu_1\left(s,t;x_0,y_0,z_0 \right) -2{\int\int_{\Gamma} K_1\left(s,t; \theta,\vartheta\right) {\mu _{1} \left(\theta,\vartheta;x_0,y_0,z_0
\right)d\theta d\vartheta}}
\end{equation*}
\begin{equation} \label{eq706}
\,\,\,\,\,\,\,\,\,\,\,\,\,\,\,\,\,\,\,\,\,\,\,\,\,\,\,\,\,\,\,\,\,\,\,\,\,\,\,\,\,\,\,\,\,
\,\,\,\,\,\,\,\,\,\,\,\,\,\,\,=2q_{1} \left(
{x(s,t),y(s,t),z(s,t);x_0,y_0,z_0} \right).
\end{equation}
Правая часть уравнения (\ref{eq706}) есть непрерывная функция от
$s$ и $t$  (точка $(x_0,y_0,z_0)$ лежит внутри $D$). В \S \ref{S5}
было доказано, что $\lambda=2$ не является собственным значением
ядра $K_1\left(s,t; \theta,\vartheta\right)$ и, следовательно,
уравнение (\ref{eq706}) разрешимо и его непрерывное решение можно
записать в виде
\begin{equation*}
 \mu_1\left(s,t;x_0,y_0,z_0 \right)=2q_{1} \left( {x(s,t),y(s,t),z(s,t);x_0,y_0,z_0} \right)+
\end{equation*}
\begin{equation} \label{eq707}
+4{\int\int_{\Gamma} R_1\left(s,t; \theta,\vartheta; 2\right)
{q_{1} \left(\xi,\eta,\zeta;x_0,y_0,z_0 \right)d\theta
d\vartheta}}
\end{equation}
где $R_1\left(s,t; \theta,\vartheta; 2\right)$ -- резольвента ядра
$K_1\left(s,t; \theta,\vartheta\right);$
$\left(x(s,t),y(s,t),z(s,t)\right)\in \Gamma$. Подставляя
(\ref{eq707}) в (\ref{eq705}), получим
\begin{equation*}
 v_1\left( x,y,z;x_0,y_0,z_0 \right) =
2{\int\int_{\Gamma}  {q_{1} \left( {\xi,\eta,\zeta;x_0,y_0,z_0}
\right)B_\nu^{\alpha}[q_{1} \left( {\xi,\eta,\zeta;x,y,z}
\right)]d\theta d\vartheta}}+
\end{equation*}
\begin{equation*}
\label{eq708} +4{\int\int_{\Gamma}\int\int_{\Gamma}
{B_\nu^{\alpha}[q_{1} \left( {\xi,\eta,\zeta;x,y,z}
\right)]R_1\left(\theta,\vartheta;s,t; 2\right)}}\times
\end{equation*}
\begin{equation}
\label{eq708} \times q_{1} \left(x(s,t),y(s,t),z(s,t);x_0,y_0,z_0
\right)d\theta d\vartheta ds dt.
\end{equation}

Определим теперь функцию
\begin{equation} \label{eq709}
g(x,y,z)= \left\{ {{\begin{array}{*{20}c}
 { v_1(x,y,z;x_0,y_0,z_0),\,\,\,\,(x,y,z) \in D,} \hfill \\
 { -q_1(x,y,z;x_0,y_0,z_0),\,\,\,\,(x,y,z) \in {D'}.} \hfill
\end{array}}}  \right.
\end{equation}

Функция $g(x,y,z)$ является регулярным решением уравнения
(\ref{eq1}) как внутри области $D$, так и внутри ${D'}$ и равна
нулю на бесконечности. Так как точка $(x_0,y_0,z_0)$ лежит внутри
$D$, то в ${D'}$ функция $g(x,y,z)$ имеет производные любого
порядка, непрерывные вплоть до $\Gamma$. Мы можем рассматривать
$g(x,y,z)$ в ${D'}$ как решение уравнения (\ref{eq1}),
удовлетворяющее граничным условиям
\begin{equation*}\label{eq710}
  \left. B_n^\alpha \left[g(x,y,z)\right]\right|_\Gamma=-B_n^\alpha\left[q_1(x(s,t),y(s,t),z(s,t);x_0,y_0,z_0)\right],
\end{equation*}

\begin{equation*}\label{eq711}
  \left.\left(x^{2\alpha}\frac{\partial g(x,y,z)}{\partial x}\right)\right|_{x=0}=0.
\end{equation*}
Это решение представим в виде потенциала простого слоя
\begin{equation}\label{eq712}
g(x,y,z)=\int\int_\Gamma
\rho_1(\theta,\vartheta;x_0,y_0,z_0)q_1(\xi,\eta,\zeta;x,y,z)d\theta
d\vartheta,\,\,\,(x,y,z)\in R_3^+\backslash \overline{D}
\end{equation}
с неизвестной плотностью $\rho_1(\theta,\vartheta;x_0,y_0,z_0)$.

Воспользовавшись второй из формул (\ref{eq404}), получим
интегральное уравнение для плотности  $\rho_1(s,t;x_0,y_0,z_0)$
\begin{equation*}
\rho_1(s,t;x_0,y_0,z_0)-2\int\int_\Gamma
K_1(\theta,\vartheta;s,t)\rho_1(\theta,\vartheta;x_0,y_0,z_0)d\theta
d\vartheta=
\end{equation*}
\begin{equation}\label{eq713}
=2B_n^\alpha\left[q_1(x(s,t),y(s,t),z(s,t);x_0,y_0,z_0)\right].
\end{equation}
Уравнение (\ref{eq713}) союзное с уравнением (\ref{eq706}). Его
правая часть есть непрерывная функция от $s$ и $t$. Таким образом,
уравнение (\ref{eq713}) имеет непрерывное решение:
\begin{equation*}
\rho_1(s,t;x_0,y_0,z_0)=2B_n^\alpha\left[q_1(x(s,t),y(s,t),z(s,t);x_0,y_0,z_0)\right]
\end{equation*}
\begin{equation}\label{eq714}
+4\int\int_\Gamma
R_1(\theta,\vartheta;s,t;2)B_\nu^\alpha\left[q_1(\xi,\eta,\zeta;x_0,y_0,z_0)\right]d\theta
d\vartheta.
\end{equation}
Значения потенциала простого слоя $g(x,y,z)$ на поверхности
$\Gamma$ равны $-q_1(x,y,z;x_0,y_0,z_0),$ т.е. такие же, как и
функции $v_1(x,y,z;x_0,y_0,z_0)$, а на плоскости $yOz$ их частные
производные по $x$ равны нулю. Отсюда в силу теоремы
единственности задачи Хольмгрена следует, что формула
(\ref{eq712}) для функции $g(x,y,z)$, определенной равенством
(\ref{eq709}), справедлива во всем полупространстве $x\geq 0$,
т.е.

$v_1(x,y,z; x_0,y_0,z_0)=$
\begin{equation}\label{eq715}
=\int\int_\Gamma
\rho_1(\theta,\vartheta;x_0,y_0,z_0)q_1(\xi,\eta,\zeta;x,y,z)d\theta
d\vartheta,\,\,\,(x,y,z)\in D.
\end{equation}
Таким образом, регулярная часть $v_1(x,y,z; x_0,y_0,z_0)$ функции
Грина представима в виде потенциала простого слоя.

Применяя первую из формул (\ref{eq404}) к (\ref{eq715}), получим
\begin{equation*}
  2B_n^{\alpha}\left[{v_1\left( x(s,t),y(s,t),z(s,t);x_0,y_0,z_0 \right)}\right]_i
\end{equation*}
\begin{equation*}
  =\rho_1(s,t;x_0,y_0,z_0)+2\int\int_\Gamma K_1(\theta,\vartheta; s,t)\rho_1(\theta,\vartheta;x_0,y_0,z_0)d\theta d\vartheta,
\end{equation*}
но, согласно (\ref{eq713}), имеем
\begin{equation*}
  2B_n^{\alpha}\left[{q_1\left( x(s,t),y(s,t),z(s,t);x_0,y_0,z_0 \right)}\right]_i
\end{equation*}

\begin{equation*}
  =\rho_1(s,t;x_0,y_0,z_0)-2\int\int_\Gamma K_1(\theta,\vartheta; s,t)\rho_1(\theta,\vartheta;x_0,y_0,z_0)d\theta d\vartheta.
\end{equation*}
Складывая почленно последние два равенства и принимая во внимание
(\ref{eq702}), будем иметь
\begin{equation}\label{eq716}
 B_n^{\alpha}\left[{G_1\left( x(s,t),y(s,t),z(s,t);x_0,y_0,z_0 \right)}\right]= \rho_1(s,t;x_0,y_0,z_0),
\end{equation}
и, следовательно, формулу (\ref{eq715}) можно записать в виде

$v_1(x,y,z; x_0,y_0,z_0)=$
\begin{equation*}\label{eq717}
=\int\int_\Gamma B_\nu^{\alpha}\left[{G_1\left(
\xi,\eta,\zeta;x_0,y_0,z_0
\right)}\right]q_1(\xi,\eta,\zeta;x,y,z)d\theta d\vartheta.
\end{equation*}

Умножая теперь обе части равенства (\ref{eq714}) на ${q_1\left(
x(s,t),y(s,t),z(s,t);x,y,z \right)}$, интегрируя по $\Gamma$  и
учитывая (\ref{eq707}) и (\ref{eq705}), получим
\begin{equation*}
v_1(x_0,y_0,z_0;x,y,z)=\int\int_\Gamma
\rho_1(\theta,\vartheta;x_0,y_0,z_0)q_1(\xi,\eta,\zeta;x,y,z)d\theta
d\vartheta.
\end{equation*}
Сравнивая это с формулой (\ref{eq715}), будем иметь
\begin{equation}\label{eq718}
 v_1(x,y,z; x_0,y_0,z_0)=v_1(x_0,y_0,z_0;x,y,z).
\end{equation}
если точки $(x,y,z)$ и $(x_0,y_0,z_0)$ находятся внутри области
$D$.

 \begin{lemma} \label{L8}{Функция Грина $G_1(x,y,z; x_0,y_0,z_0)$ симметрична от\-но\-си\-тель\-но точек $(x,y,z)$ и $(x_0,y_0,z_0)$, если они находятся внутри области $D$.}
\end{lemma}

Доказательство леммы следует из представления (\ref{eq702})
функции Грина и равенства (\ref{eq718}).

Для области $D_0$, ограниченной кругом $y^2+z^2\leq a^2$ плоскости
$yOz$ и полусферой $x^2+y^2+z^2=a^2\,(x \geq 0),$ функция Грина
задачи Хольмгрена имеет вид
\begin{equation}\label{eq719}
G_{01}(x,y,z;x_0,y_0,z_0)=q_{1}(x,y,z;x_0,y_0,z_0)-\left(\frac{a}{R}\right)^{1+2\alpha}q_{1}(x,y,z;\bar{x}_0,\bar{y}_0,\bar{z}_0),
\end{equation}
где
\begin{equation*}\label{eq720}
R^2=x_0^2+y_0^2+z_0^2,\,\,\,\bar{x}_0=\frac{a^2}{R^2}x_0,\,\,\,\bar{y}_0=\frac{a^2}{R^2}y_0,\,\,\,\bar{z}_0=\frac{a^2}{R^2}z_0.
\end{equation*}

Покажем, что регулярную часть
$$v_{01}(x,y,z;x_0,y_0,z_0)=-\left(\frac{a}{R}\right)^{1+2\alpha}q_{1}(x,y,z;\bar{x}_0,\bar{y}_0,\bar{z}_0)
$$
функции Грина $G_{01}(x,y,z;x_0,y_0,z_0)$ можно представить в виде

$v_{01}(x,y,z; x_0,y_0,z_0)=$
\begin{equation}\label{eq721}
=-\int\int_\Gamma
\rho_1(s,t;x,y,z)v_{01}(x(s,t),y(s,t),z(s,t);x_0,y_0,z_0)ds dt,
\end{equation}
где $\rho_1(s,t;x,y,z)$  есть решение уравнения (\ref{eq715}).

Действительно, пусть $(x_0,y_0,z_0)$ произвольная точка внутри
области $D$. Рассмотрим функцию
\begin{equation*}
u(x,y,z; x_0,y_0,z_0)=-\int\int_\Gamma
\rho_1(s,t;x,y,z)v_{01}(x(s,t),y(s,t),z(s,t);x_0,y_0,z_0)ds dt.
\end{equation*}
Она, как функция от $(x,y,z)$, удовлетворяет уравнению
(\ref{eq1}), так как этому уравнению удовлетворяет функция
$\rho_1(s,t;x,y,z)$. Подставляя вместо $\rho_1(s,t;x,y,z)$  ее
выражение (\ref{eq714}), получим

$u(x,y,z; x_0,y_0,z_0)=$
\begin{equation}\label{eq722}
=-\int\int_\Gamma
\psi(s,t;x_0,y_0,z_0)B_n^\alpha\left[q_{1}(x(s,t),y(s,t),z(s,t);x,y,z)\right]ds
dt,
\end{equation}
где \\
$
\psi(s,t;x_0,y_0,z_0)=2v_{01}\left(x(s,t),y(s,t),z(s,t);x_0,y_0,z_0\right)$
$$+4\int\int_\Gamma R_1(s,t;\theta,\vartheta;2)v_{01}\left(\xi,\eta,\zeta; x_0,y_0,z_0\right)d\theta d\vartheta,
$$
т.е. $\psi(s,t;x_0,y_0,z_0)$ есть решение интегрального уравнения
\begin{equation*}
  \psi(s,t;x_0,y_0,z_0)-2\int\int_\Gamma K_1(s,t;\theta,\vartheta)\psi(\theta,\vartheta; x_0,y_0,z_0)d\theta d\vartheta=
\end{equation*}
\begin{equation}\label{eq723}
 = 2v_{01}\left(x(s,t),y(s,t),z(s,t);x_0,y_0,z_0\right).
\end{equation}
Применяя первую из формул (\ref{eq324}) к потенциалу двойного слоя
(\ref{eq722}), получим
\\

$u_i\left(x(s,t),y(s,t),z(s,t); x_0,y_0,z_0\right)=$

  $$=\frac{1}{2}\psi(s,t;x_0,y_0,z_0)-\int\int_\Gamma K_1(s,t;\theta,\vartheta)\psi(\theta,\vartheta;x_0,y_0,z_0)d\theta d\vartheta,
$$
откуда в силу (\ref{eq723}) имеем
\begin{equation*}
u_i\left(x(s,t),y(s,t),z(s,t);
x_0,y_0,z_0\right)=v_{01}\left(x(s,t),y(s,t),z(s,t);
x_0,y_0,z_0\right),
\end{equation*}
\begin{equation*}
  \left(x(s,t),y(s,t),z(s,t)\right)\in \Gamma.
\end{equation*}
Нетрудно убедиться, что
\begin{equation*}
  \left.\left(x^{2\alpha}\frac{\partial u\left(x,y,z; x_0,y_0,z_0\right)}{\partial x}\right)\right|_{x=0}=0,\,\,\,\,\,\,\left.\left(x^{2\alpha}\frac{\partial v_{01}\left(x,y,z; x_0,y_0,z_0\right)}{\partial x}\right)\right|_{x=0}=0.
\end{equation*}
Таким образом, функции  $u\left(x,y,z; x_0,y_0,z_0\right)$  и
$v_{01}\left(x,y,z; x_0,y_0,z_0\right)$  удов\-ле\-твор\-яют
одному и тому же уравнению (\ref{eq1}) и одинаковым краевым
условиям, и в силу единственности решения задачи Хольмгрена
\begin{equation*}
  u\left(x,y,z; x_0,y_0,z_0\right) \equiv v_{01}\left(x,y,z; x_0,y_0,z_0\right).
\end{equation*}
Вычитая теперь из (\ref{eq702}) обе части равенства (\ref{eq719}),
получим
\begin{equation*}
 H_1\left(x,y,z; x_0,y_0,z_0\right) =  G_{1}\left(x,y,z; x_0,y_0,z_0\right) - G_{01}\left(x,y,z; x_0,y_0,z_0\right)=
\end{equation*}
\begin{equation*}
=v_{1}\left(x,y,z; x_0,y_0,z_0\right) - v_{01}\left(x,y,z;
x_0,y_0,z_0\right)
\end{equation*}
или в силу (\ref{eq715}), (\ref{eq718}), (\ref{eq719}) и
(\ref{eq721})
\begin{equation}\label{eq724}
H_1(x,y,z; x_0,y_0,z_0)=\int\int_\Gamma
\rho_1(\theta,\vartheta;x,y,z)G_{01}(\xi,\eta,\zeta;x_0,y_0,z_0)d\theta
d\vartheta.
\end{equation}

\subsection{Решение задачи Хольмгрена для уравнения (\ref{eq1})}\label{s7.2}
Пусть $(x_0,y_0,z_0)$  -- точка внутри области $D$. Рассмотрим
область $D_{\varepsilon,\delta}\subset D$, ограниченную
поверхностью $\Gamma_\varepsilon$, параллельной поверхности
$\Gamma$, и областью $X_{\varepsilon,\delta}$, лежащей на
плоскости $x=\delta>\varepsilon$. Выберем $\varepsilon$ и $\delta$
столь малыми, чтобы точка $(x_0,y_0,z_0)$ находилась внутри
$D_{\varepsilon,\delta}$. Вырежем из области
$D_{\varepsilon,\delta}$  шар малого радиуса $\rho$ с центром в
точке $(x_0,y_0,z_0)$  и оставшуюся часть $D_{\varepsilon,\delta}$
обозначим через $D_{\varepsilon,\delta}^\rho$, в которой функция
Грина $G_1(x,y,z; x_0,y_0,z_0)$  будет регулярным решением
уравнения (\ref{eq1}).

Пусть $u(x,y,z)$ есть регулярное решение уравнения (\ref{eq1}) в
области $D$, удовлетворяющее граничным условиям (\ref{eq602}).

Применяя формулу (\ref{eq9}), получим
\begin{equation*}
\int\int_{\Gamma_\varepsilon}\left(G_1 B_n^\alpha[u]-u
B_n^\alpha[G_1]\right)dsdt
+\int\int_{X_{\varepsilon,\delta}}\left.x^{2\alpha}\left(u
\frac{\partial G_1}{\partial x}-G_1\frac{\partial u}{\partial
x}\right)\right|_{x=\delta}dydz=
\end{equation*}
\begin{equation*}
  =\int\int_{C_\rho}\left(G_1 B_n^\alpha[u]-u B_n^\alpha[G_1]\right)dsdt,
\end{equation*}
где $C_\rho$ -- сфера вырезанного шара. Переходя к пределу при
$\rho \to 0$, а затем при $\varepsilon \to 0$ и $\delta \to 0$,
получим
\begin{equation*}
 u\left(x_0,y_0,z_0\right)= - \int\int_{X}\nu_1(y,z)G_1\left(0,y,z; x_0,y_0,z_0\right)dydz
 \end{equation*}
\begin{equation}\label{eq726}
 -\int\int_{\Gamma} \varphi (\theta,\vartheta) B_n^\alpha[G_1\left(\xi,\eta,\zeta;x_0,y_0,z_0\right)]d\theta d\vartheta.
\end{equation}

Покажем, что формула (\ref{eq726}) дает решение задачи Хольмгрена.

Нетрудно видеть, что первый интеграл $I_1(x_0,y_0,z_0)$ в формуле
(\ref{eq726}) есть регулярное в области $D$ решение уравнения
(\ref{eq1}), непрерывное в $\overline{D}$. Обозначим
\begin{equation*}
\varphi(x_0,y_0,z_0)=\int\int_{X}\nu_1(y,z) q_1(0,y,z;
x_0,y_0,z_0) dy dz=
\end{equation*}

\begin{equation}\label{eq727}
=\frac{1}{2\pi}\int\int_{X}\nu_1(y,z)\left[x_0^2+(y-y_0)^2+(z-z_0)^2\right]^{-\frac{1}{2}-\alpha}dy
dz.
\end{equation}
Легко видеть, что $\varphi(x_0,y_0,z_0)$ есть непрерывная функция
в $\overline{D}$. Интеграл $I_1(x_0,y_0,z_0)$ в силу (\ref{eq727})
и (\ref{eq708}) и симметричности функции $v_1(x,y,z; x_0,y_0,z_0)$
можно представить в виде

\begin{equation*}
 I_1\left(x_0,y_0,z_0 \right) =-\varphi\left(x_0,y_0,z_0 \right)-
2{\int\int_{\Gamma} \phi\left(\xi,\eta,\zeta \right)
B_\nu^{\alpha}[q_{1} \left( {\xi,\eta,\zeta;x_0,y_0,z_0}
\right)]d\theta d\vartheta}-
\end{equation*}
\begin{equation*}
\label{eq728}
-4{\int\int_{\Gamma}\int\int_{\Gamma}R_1\left(\theta,\vartheta;s,t;
2\right)\varphi\left(x(s,t),y(s,t),z(s,t) \right)}\times
\end{equation*}
\begin{equation}
\label{eq728} \times {B_\nu^{\alpha}[q_{1} \left(
{\xi,\eta,\zeta;x_0,y_0,z_0} \right)]d\theta d\vartheta ds dt}.
\end{equation}

Последние два интеграла в формуле  (\ref{eq728}) суть потенциалы
двойного слоя. Принимая во внимание первую из формул (\ref{eq324})
и интегральное уравнение для резольвенты
$R_1(s,t;\theta,\vartheta;2)$, мы из формулы (\ref{eq728}) получим
\begin{equation*} \label{eq729}
  \left.I_1\left(x_0,y_0,z_0\right)\right|_{\Gamma}=0,
\end{equation*}

Нетрудно видеть, что
\begin{equation*}
\label{eq730}   {\mathop {\lim} \limits_{x_0\to
0}}x_0^{2\alpha}\frac{\partial
I_1\left(x_0,y_0,z_0\right)}{\partial x_0} =
\nu_1\left(y_0,z_0\right),\,\,\,\,\left(y_0,z_0\right)\in X.
\end{equation*}
В самом деле, интеграл $I_1\left(x_0,y_0,z_0\right)$  в силу
(\ref{eq715}) и симметричности функции $v_1\left(x,y,z;
x_0,y_0,z_0\right)$  можно записать также в виде
\begin{equation*}
I_1(x_0,y_0,z_0)=-\int\int_{X}\nu_1(y,z) q_1(0,y,z; x_0,y_0,z_0)
dy dz-
\end{equation*}

\begin{equation*}\label{eq731}
-\int\int_{X}\nu_1(y,z)dydz\int\int_\Gamma \rho
(\theta,\vartheta;0,y,z)q_1(\xi,\eta,\zeta; x_0,y_0,z_0)d\theta
d\vartheta.
\end{equation*}
В \S \ref{S6} было показано, что производная по $x_0$ от первого
слагаемого, умноженного на $x_0^{2\alpha}$, равна $\nu_1(y_0,z_0)$
при $x_0 \to 0$, $(y_0,z_0)\in X$. Производная по $x_0$ от второго
слагаемого, умноженного на $x_0^{2\alpha}$, равна нулю при
$x_0=0$, так как $x_0^{2\alpha}\displaystyle\frac{\partial
q_1}{\partial x_0} =0$  при $x_0 \to 0$, $(y_0,z_0)\in X$.

Рассмотрим второй интеграл $I_2(x_0,y_0,z_0)$ в формуле
(\ref{eq726}), который в силу (\ref{eq716}) и (\ref{eq714}) можно
записать в виде
\begin{equation*}
 I_2(x_0,y_0,z_0)=-\int\int_{\Gamma}\varphi(s,t) \rho_1(s,t; x_0,y_0,z_0) ds dt=
\end{equation*}

\begin{equation}\label{eq732}
 =-\int\int_{\Gamma}\chi\left(\theta,\vartheta\right) B_\nu^\alpha \left[q_1(\xi,\eta,\zeta; x_0,y_0,z_0)\right] d\theta d\vartheta,
\end{equation}
где
\begin{equation*}
\chi\left(\theta,\vartheta\right)=2\varphi\left(\theta,\vartheta\right)+4\int\int_\Gamma
R_1\left(\theta,\vartheta; s,t;2\right)\varphi(s,t) ds dt,
\end{equation*}
т.е. функция $\chi\left(\theta,\vartheta\right)$  есть решение
интегрального уравнения
\begin{equation}\label{eq733}
 \chi\left(s,t\right)-2\int\int_\Gamma K_1\left(s,t;\theta,\vartheta;2\right)\chi\left(\theta,\vartheta\right) d\theta d\vartheta =2\varphi\left(s,t\right).
\end{equation}

Так как $\chi\left(s,t\right)$ -- непрерывная функция, то
$I_2(x_0,y_0,z_0)$ есть регулярное в области $D$ решение уравнения
(\ref{eq1}), непрерывное в $\overline{D}$, которое в силу
(\ref{eq324}) и (\ref{eq733}) удовлетворяет условию
\begin{equation*}
  \left.I_2\left(x_0,y_0,z_0\right)\right|_{\Gamma}=\varphi\left(s,t\right),
\end{equation*}
Нетрудно видеть, что
\begin{equation*} {\mathop {\lim} \limits_{x_0\to 0}}x_0^{2\alpha}\frac{\partial I_2\left(x_0,y_0,z_0\right)}{\partial x_0} = 0,\,\,\,\,\left(y_0,z_0\right)\in X.
\end{equation*}

Используя формулы (\ref{eq724}) и (\ref{eq719}), решение
(\ref{eq726}) задачи Хольмгрена для уравнения (\ref{eq1}) можно
записать в виде

$u\left(x_0,y_0,z_0\right)=$
\begin{equation*}
 =-\int\int_X \nu_1(y,z)\left[G_{01}\left(0,y,z; x_0,y_0,z_0\right)+H_{1}\left(0,y,z; x_0,y_0,z_0\right)\right] dy dz-
\end{equation*}
\begin{equation}\label{eq734}
 -\int\int_\Gamma \varphi(\theta,\vartheta)\left\{B_\nu^\alpha\left[G_{01}\left(\xi,\eta,\zeta; x_0,y_0,z_0\right)\right]+B_\nu^\alpha\left[H_{1}\left(\xi,\eta,\zeta; x_0,y_0,z_0\right)\right]\right\} d\theta d\vartheta,
\end{equation}
где
\begin{equation*}
G_{01}\left(0,y,z;
x_0,y_0,z_0\right)=\frac{1}{2\pi}\left\{\left[x_0^2+(y-y_0)^2+(z-z_0)^2\right]^{-\frac{1}{2}-\alpha}-\left[\left(a-\frac{yy_0}{a}\right)^2+\right.\right.
\end{equation*}

\begin{equation*}\label{eq735}
\left.\left.+\left(a-\frac{zz_0}{a}\right)^2
+\frac{1}{a^2}\left(x_0^2y^2+y^2z_0^2+x_0^2z^2+y_0^2z^2\right)-a^2\right]^{-\frac{1}{2}-\alpha}\right\}
\end{equation*}

\begin{equation*}\label{eq{736}}
H_{1}\left(x,y,z;
x_0,y_0,z_0\right)=\int\int_\Gamma\rho_1\left(\theta,\vartheta;
x_0,y_0,z_0\right)G_{01}\left(\xi,\eta,\zeta; x,y,z\right)d\theta
d\vartheta.
\end{equation*}

Решение (\ref{eq734}) задачи Хольмгрена более удобно для
дальнейших исследований. В случае полусферической области $D_0$
функция $H_{1}\left(x,y,z; x_0,y_0,z_0\right)\equiv 0$ и решение
(\ref{eq734}) принимает более простой вид:
\begin{equation*}
 u\left(x_0,y_0,z_0\right)=-\frac{1}{2\pi}\int\limits_{-a}^a dy\int\limits_{-\sqrt{a^2-y^2}}^{\sqrt{a^2-y^2}}\nu_1(y,z)\left\{\left[x_0^2+(y-y_0)^2+(z-z_0)^2\right]^{-\frac{1}{2}-\alpha}-\right.
\end{equation*}
\begin{equation*}
 -\left.\left[\left(b-\frac{yy_0}{b}\right)^2+\left(b-\frac{zz_0}{b}\right)^2
+\frac{1}{b^2}\left(x_0^2y^2+y^2z_0^2+x_0^2z^2+y_0^2z^2\right)-b^2\right]^{-\frac{1}{2}-\alpha}\right\}dz+
\end{equation*}
\begin{equation}\label{eq737}
+\frac{1+2\alpha}{2\pi}\int\int_\Gamma
\varphi\left(\theta,\vartheta\right)
\xi^{2\alpha}F\left(\frac{3}{2}+\alpha,\alpha;2\alpha;1-\frac{r_1^2}{r^2}\right)\frac{c^2-R^2}{cr^{3+2\alpha}}d\theta
d\vartheta,
\end{equation}
где
\begin{equation*}
0<2\alpha<1; \,\,b^2=y^2+z^2,\,\,\,c^2=\xi^2+\eta^2+\zeta^2,
\,\,\xi>0;\,\,\, R^2=x_0^2+y_0^2+z_0^2, \,\,x_0>0;
\end{equation*}
\begin{equation*}
r^2=\left(\xi-x_0\right)^2+\left(\eta-y_0\right)^2+\left(\zeta-z_0\right)^2;
\,\,r_1^2=\left(\xi+x_0\right)^2+\left(\eta-y_0\right)^2+\left(\zeta-z_0\right)^2;
\end{equation*}
\begin{equation*}
\xi=\xi(\theta,\vartheta),\,\,\eta=\eta(\theta,\vartheta),\,\,\zeta=\zeta(\theta,\vartheta),\,\,(\xi,\eta,\zeta)\in\Gamma.
\end{equation*}

Непосредственным вычислением можно показать, что функция,
определенная формулой  (\ref{eq734}) является решением задачи
Хольмгрена для уравнения (\ref{eq1}), т.е. она удовлетворяет
уравнению (\ref{eq1}) и условиям (\ref{eq602}).

\subsection{Функция Грина задачи Дирихле}\label{s7.3}

\begin{definition}\label{de4}
Функцией Грина задачи Дирихле для уравнения (\ref{eq1}) называется
функция $G_2(x,y,z;x_0,y_0,z_0)$, удовлетворяющая следующим
условиям:

1) внутри области $D$, кроме точки $(x_0,y_0,z_0)$, эта функция
есть регулярное решение уравнения (\ref{eq1});

2) она удовлетворяет граничному  условию
\begin{equation}\label{eq738}
  \left. G_2(x,y,z;x_0,y_0,z_0)\right|_{\Gamma\cup \bar{X}}=0;
\end{equation}

3) она может быть представлена в виде
\begin{equation}\label{eq739}
G_2(x,y,z;x_0,y_0,z_0)=q_2(x,y,z;x_0,y_0,z_0)+v_2(x,y,z;x_0,y_0,z_0)
\end{equation}
где $q_{2} \left( {x,y,z;\xi,\eta,\zeta}  \right)$ --
фундаментальное решение уравнения (\ref{eq1}), определенное
формулой (\ref{eq3}), а $v_2(x,y,z;x_0,y_0,z_0)$  -- регулярное
решение уравнения (\ref{eq1}) везде внутри $D$. \end{definition}

Построение функции Грина сводится к нахождению ее регулярной части
$v_2(x,y,z;x_0,y_0,z_0)$, которая в силу (\ref{eq6}),\,
(\ref{eq738}) и (\ref{eq739}) должна удовлетворять граничным
условиям
\begin{equation*}\label{eq741}
  \left. v_2(x,y,z;x_0,y_0,z_0)\right|_\Gamma=\left.-q_2(x,y,z;x_0,y_0,z_0)\right|_\Gamma,
\end{equation*}
\begin{equation*}\label{eq742}
   v_2(0,y,z;x_0,y_0,z_0)=0.
\end{equation*}

Функцию $v_2(x,y,z;x_0,y_0,z_0)$  будем искать в виде потенциала
двойного слоя:
\begin{equation}
\label{eq743} v_2\left( x,y,z;x_0,y_0,z_0 \right) =
{\int\int_{\Gamma}  {\mu _{2} \left( \theta,\vartheta;x_0,y_0,z_0
\right)B_\nu^{\alpha}[q_{2} \left( {\xi,\eta,\zeta;x,y,z}
\right)]d\theta d\vartheta}}.
\end{equation}
Воспользовавшись первой из формул (\ref{eq332}), получим
интегральное уравнение для плотности $\mu_2\left(s,t;x_0,y_0,z_0
\right)$
\begin{equation*}
 \mu_2\left(s,t;x_0,y_0,z_0 \right) -2{\int\int_{\Gamma} K_2\left(s,t; \theta,\vartheta\right) {\mu _{2} \left(\theta,\vartheta;x_0,y_0,z_0
\right)d\theta d\vartheta}}=
\end{equation*}
\begin{equation} \label{eq744}
\,\,\,\,\,\,\,\,\,\,\,\,\,\,\,\,\,\,\,\,\,\,\,\,\,\,\,\,\,\,\,\,\,\,\,\,\,\,\,\,\,\,\,\,\,
\,\,\,\,\,\,\,\,\,\,\,\,\,\,\,=2q_{2} \left(
{x(s,t),y(s,t),z(s,t);x_0,y_0,z_0} \right).
\end{equation}

Правая часть уравнения (\ref{eq744}) есть непрерывная функция от
$s$ и $t$  (точка $(x_0,y_0,z_0)$ лежит внутри $D$). В силу леммы
\ref{L7} к уравнению (\ref{eq744}) применима теория Фредгольма. В
разделе \ref{S5} было доказано, что $\lambda=2$ не является
собственным значением ядра $K_2\left(s,t; \theta,\vartheta\right)$
и, следовательно, уравнение (\ref{eq744}) разрешимо и его
непрерывное решение можно записать в виде
\begin{equation*}
 \mu_2\left(s,t;x_0,y_0,z_0 \right)=2q_{2} \left( {x(s,t),y(s,t),z(s,t);x_0,y_0,z_0} \right)
\end{equation*}
\begin{equation} \label{eq745}
+4{\int\int_{\Gamma} R_2\left(s,t; \theta,\vartheta; 2\right)
{q_{2} \left(\xi,\eta,\zeta;x_0,y_0,z_0 \right)d\theta
d\vartheta}},
\end{equation}
где $R_2\left(s,t; \theta,\vartheta; 2\right)$ -- резольвента ядра
$K_2\left(s,t; \theta,\vartheta\right);$
$\left(x(s,t),y(s,t),z(s,t)\right)\in \Gamma$. Подставляя
(\ref{eq745}) в (\ref{eq743}), получим
\begin{equation*}
 v_2\left( x,y,z;x_0,y_0,z_0 \right) =
2{\int\int_{\Gamma}  {q_{2} \left( {\xi,\eta,\zeta;x_0,y_0,z_0}
\right)B_\nu^{\alpha}[q_{2} \left( {\xi,\eta,\zeta;x,y,z}
\right)]d\theta d\vartheta}}
\end{equation*}
\begin{equation*}
\label{eq746} +4\int\int_{\Gamma}\int\int_{\Gamma}
B_\nu^{\alpha}[q_{2} \left( {\xi,\eta,\zeta;x,y,z}
\right)]R_2\left(\theta,\vartheta;s,t; 2\right)\times
\end{equation*}
\begin{equation*}
\label{eq746} \times q_{2} \left(x(s,t),y(s,t),z(s,t);x_0,y_0,z_0
\right)d\theta d\vartheta ds dt.
\end{equation*}

Совершенно так же, как и при построении функции Грина задачи
Хольмгрена, можно показать, что регулярная часть $v_2\left(
x,y,z;x_0,y_0,z_0 \right)$ функции Грина задачи Дирихле
представима в виде потенциала простого слоя
\begin{equation}\label{eq747}
v_2(x,y,z; x_0,y_0,z_0)=\int\int_\Gamma
\rho_2(\theta,\vartheta;x_0,y_0,z_0)q_2(\xi,\eta,\zeta;x,y,z)d\theta
d\vartheta,
\end{equation}
где
\begin{equation*}
\rho_2(s,t;x_0,y_0,z_0)=2B_n^\alpha\left[q_2(x(s,t),y(s,t),z(s,t);x_0,y_0,z_0)\right]
\end{equation*}
\begin{equation*}\label{eq748}
+4\int\int_\Gamma
R_2(\theta,\vartheta;s,t;2)B_n^\alpha\left[q_2(\xi,\eta,\zeta;x_0,y_0,z_0)\right]d\theta
d\vartheta,
\end{equation*}
т.е.  $\rho_2(s,t;x_0,y_0,z_0)$ есть решение интегрального
уравнения
\begin{equation*}
\rho_2(s,t;x_0,y_0,z_0)-2\int\int_\Gamma
K_2(\theta,\vartheta;s,t)\rho_2(\theta,\vartheta;x_0,y_0,z_0)d\theta
d\vartheta=
\end{equation*}
\begin{equation}\label{eq749}
=2B_n^\alpha\left[q_2(x(s,t),y(s,t),z(s,t);x_0,y_0,z_0)\right].
\end{equation}

Применяя первую из формул (\ref{eq416}) к (\ref{eq747}), получим
\begin{equation*}
  2B_n^{\alpha}\left[{v_2\left( x(s,t),y(s,t),z(s,t);x_0,y_0,z_0 \right)}\right]_i
\end{equation*}
\begin{equation*} \label{eq750}
  =\rho_2(s,t;x_0,y_0,z_0)+2\int\int_\Gamma K_2(\theta,\vartheta; s,t)\rho_1(\theta,\vartheta;x_0,y_0,z_0)d\theta d\vartheta.
\end{equation*}
Отсюда, принимая во внимание  (\ref{eq739}) и (\ref{eq749}), будем
иметь
\begin{equation*}\label{eq751}
 B_n^{\alpha}\left[{G_2\left( x(s,t),y(s,t),z(s,t);x_0,y_0,z_0 \right)}\right]= \rho_2(s,t;x_0,y_0,z_0),
\end{equation*}
и, следовательно, формулу (\ref{eq747}) можно записать в виде

$v_2(x,y,z; x_0,y_0,z_0)=$
\begin{equation*}\label{eq752}
=\int\int_\Gamma q_2(\xi,\eta,\zeta;x,y,z)
B_\nu^{\alpha}\left[{G_2\left( \xi,\eta,\zeta;x_0,y_0,z_0
\right)}\right]d\theta d\vartheta.
\end{equation*}

 \begin{lemma} \label{L9}{Функция Грина $G_2(x,y,z; x_0,y_0,z_0)$  симметрична от\-но\-си\-тель\-но точек $(x,y,z)$ и $(x_0,y_0,z_0)$, если они находятся внутри области $D$.}
\end{lemma}

Доказательство этой  леммы проводится аналогично доказательству
леммы \ref{L8}.

\subsection{Решение задачи Дирихле для уравнения (\ref{eq1})}\label{s7.4}
Для полушаровой области $D_0$, ограниченной кругом $y^2+z^2\leq
a^2$ плоскости $yOz$ и полусферой $x^2+y^2+z^2=a^2\,(x \geq 0),$
функция Грина задачи Дирихле имеет вид

$G_{02}(x,y,z;x_0,y_0,z_0)=$
\begin{equation}\label{eq753}
=q_{2}(x,y,z;x_0,y_0,z_0)-\left(\frac{a}{R}\right)^{3-2\alpha}q_{2}(x,y,z;\bar{x}_0,\bar{y}_0,\bar{z}_0),
\end{equation}
где
\begin{equation*}
a^2=x^2+y^2+z^2,\,\,\,R^2=x_0^2+y_0^2+z_0^2,\,\,\,\bar{x}_0=\frac{a^2}{R^2}x_0,\,\,\,\bar{y}_0=\frac{a^2}{R^2}y_0,\,\,\,\bar{z}_0=\frac{a^2}{R^2}z_0.
\end{equation*}

Регулярную часть
$$v_{02}(x,y,z;x_0,y_0,z_0)=-\left(\frac{a}{R}\right)^{3-2\alpha}q_{2}(x,y,z;\bar{x}_0,\bar{y}_0,\bar{z}_0)
$$
функции Грина $G_{02}(x,y,z;x_0,y_0,z_0)$ также можно представить
в виде

$v_{02}(x,y,z; x_0,y_0,z_0)=$
\begin{equation}\label{eq754}
=-\int\int_\Gamma
\rho_2(s,t;x,y,z)v_{02}(x(s,t),y(s,t),z(s,t);x_0,y_0,z_0)ds dt,
\end{equation}
где $\rho_2(s,t;x,y,z)$  есть решение уравнения (\ref{eq749}).
Вычитая теперь из (\ref{eq739}) обе части равенства (\ref{eq753})
и учитывая (\ref{eq747}), (\ref{eq754})  и (\ref{eq753}), а также
симметричность функции Грина, получим
\begin{equation*}
 H_2\left(x,y,z; x_0,y_0,z_0\right) =  G_{2}\left(x,y,z; x_0,y_0,z_0\right) - G_{02}\left(x,y,z; x_0,y_0,z_0\right)=
\end{equation*}

\begin{equation}\label{eq755}
=\int\int_\Gamma
\rho_2(\theta,\vartheta;x,y,z)G_{02}(\xi,\eta,\zeta;x_0,y_0,z_0)d\theta
d\vartheta.
\end{equation}

\begin{theorem}\label{T6} Функция
\begin{equation*}
  u\left(x_0,y_0,z_0\right)=\int\int_X\tau_1(y,z)\left.\left(x^{2\alpha}\frac{\partial G_2\left(x,y,z; x_0,y_0,z_0\right)}{\partial x}\right)\right|_{x=0}dy dz-
\end{equation*}
\begin{equation}\label{eq756}
  -\int\int_\Gamma \varphi(\theta,\vartheta)B_\nu^\alpha\left[G_2\left(\xi,\eta,\zeta; x_0,y_0,z_0\right)\right]d\theta d\vartheta,
\end{equation}
где $\tau_1(y,z)$ -- непрерывная функция при  $(y,z)\in
\overline{\Gamma},$ а $\varphi(\theta,\vartheta)$   -- непрерывная
функция при $(\theta, \vartheta)\in \overline{\Phi}$, причем
$\left.\varphi\right|_\gamma=\left.\tau_1\right|_\gamma$, есть
решение задачи Дирихле для уравнения (\ref{eq1}) в области $D$.
\end{theorem}

Эта теорема доказывается рассуждениями аналогичными тем, которые
были приведены при решении задачи Хольмгрена.

Используя формулы (\ref{eq753}) и (\ref{eq755}), решение
(\ref{eq756}) можно представить в виде
\begin{equation*}
  u\left(x_0,y_0,z_0\right)=\int\int_X\tau_1(y,z)\left[\left.\left(x^{2\alpha}\frac{\partial G_{02}}{\partial x}\right)\right|_{x=0}
  +\left.\left(x^{2\alpha}\frac{\partial H_2}{\partial x}\right)\right|_{x=0}\right]dy dz-
\end{equation*}
\begin{equation}\label{eq757}
  -\int\int_\Gamma \varphi(\theta,\vartheta)\left\{B_\nu^\alpha\left[G_{02}\left(\xi,\eta,\zeta; x_0,y_0,z_0\right)\right]+
  B_\nu^\alpha\left[H_2\left(\xi,\eta,\zeta; x_0,y_0,z_0\right)\right]\right\}d\theta d\vartheta,
\end{equation}
где
\begin{equation*}
\left.\left(x^{2\alpha}\frac{\partial G_{02}}{\partial
x}\right)\right|_{x=0}=\frac{1-2\alpha}{2\pi}
x_0^{1-2\alpha}\displaystyle\left\{\left[x_0^2+(y-y_0)^2+(z-z_0)^2\right]^{\alpha-\frac{3}{2}}-\right.
\end{equation*}
\begin{equation*}
-\left.\left[\left(a-\frac{yy_0}{a}\right)^2+\left(a-\frac{zz_0}{a}\right)^2
+\frac{1}{a^2}\left(x_0^2y^2+y^2z_0^2+x_0^2z^2+y_0^2z^2\right)-a^2\right]^{\alpha-\frac{3}{2}}\right\}.
\end{equation*}

Отметим, что решение (\ref{eq757}) более удобно для дальнейших
исследований. В случае полусферической области $D_0$ функция
$H_2\left(\xi,\eta,\zeta; x_0,y_0,z_0\right)\equiv 0$ и решение
(\ref{eq757}) принимает более простой вид:

$u\left(x_0,y_0,z_0\right)=$
\begin{equation*}
 =\frac{1-2\alpha}{2\pi} x_0^{1-2\alpha}\int\limits_{-a}^a dy\int\limits_{-\sqrt{a^2-y^2}}^{\sqrt{a^2-y^2}} \tau_1(y,z)\left\{\left[x_0^2+(y-y_0)^2+(z-z_0)^2\right]^{\alpha-\frac{3}{2}}-\right.
\end{equation*}
\begin{equation*}
-\left.\left[\left(b-\frac{yy_0}{b}\right)^2+\left(b-\frac{zz_0}{b}\right)^2
+\frac{1}{b^2}\left(x_0^2y^2+y^2z_0^2+x_0^2z^2+y_0^2z^2\right)-b^2\right]^{\alpha-\frac{3}{2}}\right\}dz+
\end{equation*}
\begin{equation}\label{eq758}
+\frac{3-2\alpha}{2\pi} x_0^{1-2\alpha}\int\int_\Gamma
\varphi\left(\theta,\vartheta\right) \xi
F\left(\frac{5}{2}-\alpha,1-\alpha;2-2\alpha;1-\frac{r_1^2}{r^2}\right)\frac{c^2-R^2}{cr^{5-2\alpha}}d\theta
d\vartheta,
\end{equation}
где
\begin{equation*}
0<2\alpha<1;\,\, b^2=y^2+z^2\,\,\,c^2=\xi^2+\eta^2+\zeta^2,
\,\,\xi>0;\,\,\, R^2=x_0^2+y_0^2+z_0^2, \,\,x_0>0;
\end{equation*}
\begin{equation*}
r^2=\left(\xi-x_0\right)^2+\left(\eta-y_0\right)^2+\left(\zeta-z_0\right)^2;
\,\,r_1^2=\left(\xi+x_0\right)^2+\left(\eta-y_0\right)^2+\left(\zeta-z_0\right)^2;
\end{equation*}
\begin{equation*}
\xi=\xi(\theta,\vartheta),\,\,\eta=\eta(\theta,\vartheta),\,\,\zeta=\zeta(\theta,\vartheta),\,\,(\xi,\eta,\zeta)\in\Gamma.
\end{equation*}

Непосредственным вычислением можно показать, что функция,
определенная формулой  (\ref{eq757}) является решением задачи
Дирихле для уравнения (\ref{eq1}), т.е. она удовлетворяет
уравнению (\ref{eq1}) и условиям (\ref{eq601}).

В заключении отметим, что полученные формулы
(\ref{eq734}),(\ref{eq737}),(\ref{eq757}) и (\ref{eq758}) играют
важную роль при изучении краевых задач для уравнения смешанного
типа, т.е. для уравнений элллиптико-гиперболического или
эллиптико-параболического типов: каждая из этих формул позволяет
легко вывести основное функциональное соотношение между следами
искомого решения и его производной на линии вырождения,
принесенное из эллиптической части смешанной области.

\end{document}